\documentclass{article}

\usepackage[a4paper, total={6in, 8in}]{geometry}
\usepackage{authblk}
\usepackage{amsmath,amssymb,amsfonts,amsthm, bm}
\usepackage{euscript}
\usepackage[utf8]{inputenc}  
\usepackage[T1]{fontenc}
\usepackage[T1]{fontenc}
\usepackage{ifpdf}
\usepackage[english]{babel}
\usepackage{caption}
\usepackage{subfig, color}
\usepackage{graphicx}
\usepackage{amsmath}
\usepackage{mathrsfs}
\usepackage{epic}
\usepackage{graphpap}
\usepackage{mathrsfs}
\usepackage{psfrag}
\usepackage{cases}
\usepackage{booktabs}
\usepackage[table]{xcolor}

\author[$\dagger$]{Stefan Bilbao}
\author[$\ddagger$]{Michele Ducceschi}
\author[$\star$]{Fabiana Zama}

\title{Explicit Exactly Energy-conserving Methods for Hamiltonian Systems}

\affil[$\dagger$]{\footnotesize{Acoustics and Audio Group,
   University of Edinburgh, Edinburgh, UK}}
 \affil[$\ddagger$]{Department of Industrial Engineering, University of Bologna, Bologna, Italy}
 \affil[$\star$]{Department of Mathematics, University of Bologna, Bologna, Italy}
\begin{document}

\maketitle

\begin{abstract}

For Hamiltonian systems, simulation algorithms that exactly conserve numerical energy or pseudo-energy have seen extensive investigation. Most available methods either require the iterative solution of nonlinear algebraic equations at each time step, or are explicit, but where the exact conservation property depends on the exact evaluation of an integral in continuous time. Under further restrictions, namely that the potential energy contribution to the Hamiltonian is non-negative, newer techniques based on invariant energy quadratisation allow for exact numerical energy conservation and yield linearly implicit updates, requiring only the solution of a linear system at each time step. In this article, it is shown that, for a general class of Hamiltonian systems, and under the non-negativity condition on potential energy, it is possible to arrive at a fully explicit method that exactly conserves numerical energy.  Furthermore, such methods are unconditionally stable, and are of comparable computational cost to the very simplest integration methods (such as St{\"o}rmer-Verlet). A variant of this scheme leading to a conditionally-stable method is also presented, and follows from a splitting of the potential energy. Various numerical results are presented, in the case of the classic test problem of Fermi, Pasta and Ulam, as well as for nonlinear systems of partial differential equations, including those describing high amplitude vibration of strings and plates. 

\end{abstract}



\section{Introduction}

The design of exact energy conserving numerical methods for nonlinear Hamiltonian systems goes back as least as far as the work of LaBudde and Greenspan, for the single particle subject to a central potential \cite{LaBudde}, Marciniak \cite{Marciniak} for the $N$-body problem, Strauss and Vazquez \cite{Strauss78} for the nonlinear Klein-Gordon equation, and Greenspan on the nonlinear harmonic oscillator \cite{Greenspan}, and was greatly generalized by Simo et al. \cite{Simo92} in their work on energy- and momentum-conserving algorithms. The concept of energy conservation in time stepping schemes goes back much further in the case of linear systems---see the discussion and references in \cite{Simo92}, particularly with regard to the energy conservation of the Crank-Nicholson method for linear systems. An exact numerical conservation of a ``pseudo-energy'' (henceforth simply ``energy'' or ``numerical energy'' in this article) of a numerical scheme for a conservative system may be contrasted with schemes that conserve energy approximately---see, e.g., early work by Hughes et. al \cite{Hughes78}. The same point of view is taken with other geometrical numerical integration techniques, such as, e.g., symplectic methods, where energy is conserved approximately, following the result from Zhong and Marsden \cite{Zhong88} that, under some restrictions, exact energy conservation and symplecticity for Hamiltonian systems cannot be obtained simultaneously \cite{Hairer06}. Various studies examine the degree to which energy is conserved for symplectic and non-symplectic methods \cite{McLachlan04, Faou}. In this paper, we focus on exact conservation of numerical energy rather than approximate conservation, and we will not consider momentum conservation. Part of the reason for the focus here on exact numerical energy conservation is the possibility of determining conditions for numerical stability as has been pointed out earlier---see \cite{Simo92} and Remark 8 in \cite{Marazzato}. 

Most exact numerical energy-conserving algorithms presented to date are implicit \cite{Strauss78, Marciniak, Simo92, vuquoc93}, requiring the solution of a nonlinear system of algebraic equations at each time step, normally through an iterative method (such as, e.g., Newton-Raphson). This is obviously a computational bottleneck, and a fully explicit formulation is preferable. Note, however, that for linear systems, explicit and exact energy-conserving methods are available---see, e.g., \cite{Diaz}. Exact explicit and conservative intergrators have also been designed on a case-by-case basis---normally new ad hoc stability considerations appear in a problem-dependent way regarding the choice of the time step. See \cite{Chin89} for a treatment of the three-body problem, and \cite{Shadwick98} for the Kepler problem. 

In a recent article by Marazzato et al. \cite{Marazzato}, an explicit Hamiltonian integrator is presented that preserves a numerical energy exactly, for general nonlinear systems. Two features are worth noting: a) the conserved quantity is defined through a continuous time integration, and thus must itself be approximated, leading effectively to approximate (non-exact) energy conservation in a discrete time implementation---with fine enough approximation, numerical energy conservation to machine precision may be achieved; and b) when the energy of the model system is non-negative, as it commonly is in practice, the numerical energy does not inherit this property, and thus conditions for numerical stability do not immediately follow. Other exactly conservative methods, also reliant on the exact evaluation of a continuous integration have been proposed previously \cite{Quispel, Brugnano}.

In the context of gradient flows for diffusive systems, recent progress has been made in developing energy-stable methods through a variable transformation representing the energy, generally a nonlinear function of the state, as a quadratic form---such methods are referred to as invariant energy quadratisation (IEQ) approaches \cite{Yang16, Yang17}. The main result is that it is possible to arrive at updating equations for a time-stepping method that are linearly implicit---so that the update depends only on the solution of a linear system, rather than the solution of a system of nonlinear algebraic equations. This allows for the sidestepping of the many difficulties associated with iterative solvers, and reduced computational cost. Alongside the more recently introduced scalar auxiliary variable (SAV) approaches \cite{SHEN2018407}, such techniques have been applied to a wide variety of problems \cite{Gong19}. More recently, IEQ/SAV approaches have been applied to Hamiltonian systems, as in the case of diagonally-implicit Runge Kutta methods \cite{Zhang20}. For an interesting overview of the relationship between quadratisation techniques and linearly implicit schemes, see the recent article by Sato et al. \cite{Sato21}.

A linearly implicit method will require the solution of a linear system at each time step, and will most likely constitute the computational bottleneck in a simulation as a whole. As will be shown here, using IEQ and SAV approaches, it is possible to arrive at energy-stable (indeed exactly lossless in machine arithmetic) numerical designs that are fully explicit. The explicit character of the update derives from the availability of a closed form inverse for the linearly implicit system that arises, and furthermore the ability to solve an $N\times N$ linear system in $O(N)$ operations, through the Sherman-Morrison inversion theorem \cite{Sherman50}. This property is fully general, and not dependent on the particular form of the Hamiltonian, except through the additional non-negativity requirement on the potential energy. As a result, computer execution time for exactly energy conserving methods for Hamiltonian systems is very nearly on par with that of the simplest explicit numerical schemes, such as, e.g., St{\"o}rmer-Verlet integration. This design addresses the two points a) and b) above with reference to the scheme presented in \cite{Marazzato}. 

In Section \ref{Hamsec}, a restricted class of Hamiltonian systems is introduced. Under the constraint that the potential energy is non-negative, quadratisation techniques can be used to arrive at a simplified system of three equations in momentum, position, and a single scalar auxiliary variable derived from the potential energy. A further generalised form is also discussed, where a quadratic term is extracted or split from the expression for the potential energy before quadratisation, leading to a distinct starting point for numerical designs. Time stepping methods are introduced in Section \ref{numsec}, including first the rudimentary St{\"o}rmer-Verlet method, the energy-conserving method presented in \cite{Marazzato}, and then an exact energy-conserving and unconditionally stable method obtained using IES/SAV approaches, and allowing an explicit update. A variation of this last scheme to the case of a split potential energy is also exactly energy-conserving, and leads to a distinct conditionally-stable method. In Section \ref{examplesec}, three examples are presented: the Fermi-Pasta-Ulam ODE system, and two PDE systems then reduced, by semi-discretisation, to Hamiltonian ODE systems: the nonlinear vibration of a string, and the high-amplitude vibration of a thin plate. Various numerical results are presented, illustrating convergence rates, exact numerical energy conservation to machine precision, and relative computation times. Some concluding remarks appear in Section \ref{concsec}. 

Preliminary results have been presented recently at the 2022 European Nonlinear Dynamics Conference \cite{bilbaoENOC2022}.

\section{Hamiltonian Systems}
\label{Hamsec}
Consider a system of particles, with $N$ generalised coordinates ${\bf q} = [q_{1},\hdots,q_{N}]^{\intercal}$ and momenta ${\bf p} = [p_{1},\hdots,p_{N}]^{\intercal}$. Here, $^{\intercal}$ indicates a transposition operation, so ${\bf p}$ and ${\bf q}$ are $N\times 1$ column vectors.  Both are functions of time $t\geq 0$, so ${\bf p} = {\bf p}(t)$ and ${\bf q} = {\bf q}(t)$. Suppose also that an associated Hamiltonian $H({\bf p},{\bf q})$ is defined by
\begin{equation}
\label{Hamdef}
    H\left({\bf p},{\bf q}\right) = \tfrac{1}{2}{\bf p}^{\intercal}{\bf M}^{-1}{\bf p}+{V}\left({\bf q}\right)
\end{equation}
for some constant symmetric positive definite $N\times N$ matrix ${\bf M}$ referred to as the mass matrix. $V({\bf q})$ is the potential energy for the system of particles. $H$ and $V$ are scalar functions of time $t$, through their dependence on ${\bf p}$ and ${\bf q}$. 

The system dynamics follow from Hamilton's equations:
\begin{equation}
    \dot{{\bf q}} = \nabla_{{\bf p}}H\qquad\qquad \dot{{\bf p}} = -\nabla_{{\bf q}}H\, ,
\end{equation}
where dots indicate differentiation with respect to time $t$, and $\nabla_{{\bf p}}$ and $\nabla_{{\bf q}}$ are gradients with respect to ${\bf p}$ and ${\bf q}$, respectively. For the particular form of Hamiltonian given in \eqref{Hamdef}, Hamilton's equations become
\begin{equation}
\label{Hamdyn}
     \dot{{\bf q}} = {\bf M}^{-1}{\bf p}\qquad\qquad\dot{{\bf p}} = -\nabla_{{\bf q}}V\, .
\end{equation}
Through time differentiation of the first of \eqref{Hamdyn}, and substitution of the second, a second order system of ordinary differential equations in ${\bf q}$ results:
\begin{equation}
\label{Ham_second}
    {\bf M}\ddot{{\bf q}} = -\nabla_{{\bf q}}V\, .
\end{equation}
Such a second order form serves as the starting point for many numerical methods, including the classic St{\"o}rmer-Verlet method. See Section \ref{Stormersec}. In the first order representation \eqref{Hamdyn}, initial conditions are required for ${\bf q}$ and ${\bf p}$, so ${\bf q}(0) = {\bf q}^{(0)}$ and ${\bf p}(0) = {\bf p}^{(0)}$, for given $N$ vectors ${\bf q}^{(0)}$ and ${\bf p}^{(0)}$. For the second order form \eqref{Ham_second}, one may set $\dot{{\bf q}}(0) = {\bf M}^{-1}{\bf p}^{(0)}$. 

The defining feature of a Hamiltonian system is that, given initial conditions ${\bf p}^{(0)}$ and ${\bf q}^{(0)}$, the energy $H(t)$ is constant for all $t\geq 0$:
\begin{equation}
\label{Hamconst}
    \dot{H} = 0\quad\rightarrow\qquad H(t) = H(0) = {\rm constant}\quad\forall t\geq 0 \, ,
\end{equation}
where $H(0) = H({\bf p}^{(0)}, {\bf q}^{(0)})$, evaluated from \eqref{Hamdef} using initial conditions ${\bf q}^{(0)}$ and ${\bf p}^{(0)}$. 

\subsection{Comments}

The system and Hamiltonian in \eqref{Hamdef} are not the most general possible. A few remarks are offered here:
\begin{itemize}
    \item The mass matrix ${\bf M}$ is assumed constant---a common assumption \cite{Simo92}. In many cases of interest, it is furthermore diagonal, but we will consider the more general form here, and indicate cases in which a diagonal form of ${\bf M}$ will have an impact on computational performance. 
    \item The Hamiltonian in \eqref{Hamdef} is separable, so that $H({\bf p}, {\bf q}) = T({\bf p})+V({\bf q})$, for kinetic energy $T$ and potential energy $V$. In some cases of interest \cite{Hairer06}, the mass matrix ${\bf M}$ may be of the form ${\bf M} = {\bf M}({\bf q})$, disturbing this splitting, but these cases will not be considered here.
    \item Following from the points above, the Hamiltonian is quadratic in ${\bf p}$, meaning that the associated dynamical system is linear in ${\bf p}$. 
    \item The individual displacements $q_{i}$ and momenta $p_{i}$, $i=1,\hdots,N$, are assumed scalar here, and represent either displacements/momenta in a signal coordinate direction, or individual components of more general vector displacements/momenta. Both cases will be seen in the numerical examples in Section \ref{examplesec}. 
    
\end{itemize}

\subsection{Non-negativity of the Potential Energy}

As a further constraint, assume that
\begin{equation}
\label{Vpos}
    V({\bf q})\geq 0\qquad \forall {\bf q}\in{\mathbb R}^{N}\, .
\end{equation}
This constraint is a natural one in many applications, but not all---for example, the gravitational potential used in the calculations of planetary orbits is not of this form \cite{Gonzalez96}. This further implies, from \eqref{Hamdef}, that 
\begin{equation}
    H\geq 0\qquad \forall {\bf p}, {\bf q}\in {\mathbb R}^{N}\, .
\end{equation}
Furthermore, as ${\bf M}>0$, from \eqref{Hamconst}, one has the following bound on ${\bf p}(t)$ in terms of the initial conditions:
\begin{equation}
\label{pbound}
    \|{\bf p}(t)\|\leq \sqrt{2\lambda_{{\rm max}}({\bf M})H(0)}\qquad \forall t\geq 0\, ,
\end{equation}
where $\lambda_{{\rm max}}({\bf M})$ is the maximal eigenvalue of ${\bf M}$, and $\|\cdot\|$ indicates a Euclidean norm. If $V({\bf q})$ is radially unbounded \cite{terrell}, so that $V({\bf q})\rightarrow +\infty$ as $\|{\bf q}\|\rightarrow +\infty$, then a further bound follows for $\|{\bf q}\|$. 

The non-negativity property of $V$ is essential to invariant energy quadratisation methods---indeed, it can be generalized to the case of $V$ bounded from below, so that $V({\bf q})\geq c$, for any constant $c$ \cite{SHEN2018407}, but the non-negativity condition \eqref{Vpos} above will be used here for simplicity. 

\subsection{Potential Energy Quadratisation}

Under the non-negativity condition on $V$, from \eqref{Vpos}, one may define
\begin{equation}
\label{psidef}
    V = \tfrac{1}{2}\psi^2\, .
\end{equation}
Hamilton's equations become
\begin{equation}
\label{Hamdyn_quad}
    \dot{{\bf q}} = {\bf M}^{-1}{\bf p}\qquad\qquad \dot{{\bf p}} = -\psi\nabla_{{\bf q}}\psi\, .
\end{equation}
Furthermore, using the chain rule, 
\begin{equation}
    \dot{\psi} = (\nabla_{{\bf q}}\psi)^{\intercal}\dot{{\bf q}}\, .
\end{equation}
Finally, introducing 
\begin{equation}
\label{gdef}
{\bf g} \triangleq \nabla_{{\bf q}}\psi = \frac{1}{\sqrt{2V}}\nabla_{{\bf q}}V\, ,
\end{equation} 
a system of three equations results: 
\begin{equation}
\label{Hamdyn_g}
    \dot{{\bf q}} = {\bf M}^{-1}{\bf p}\qquad\qquad \dot{{\bf p}} = -\psi{\bf g}\qquad\qquad \dot{\psi} = {\bf g}^{\intercal}\dot{{\bf q}}\, .
\end{equation}
An auxiliary initial condition $\psi(0)$ may be set as $\psi(0) = \sqrt{2V({\bf q}^{(0)})}$ in terms of the given initial condition ${\bf q}(0) = {\bf q}^{(0)}$. 

The key feature of system \eqref{Hamdyn_g} is that, if ${\bf g}$ is assumed known at at any given time instant, the resulting equations are linear in ${\bf p}$, ${\bf q}$ and $\psi$. Though in the present case, ${\bf g}$ does indeed have a dependence on ${\bf q}$, in the numerical setting, ${\bf g}$ can be evaluated directly using previously computed values of the solution, and thus the update becomes a matter of solving a linear system. This is the essence of the ``linearly implicit'' property of schemes arising from energy quadratisation. As will be shown subsequently here, under the appropriate numerical design, the linear system is of a particularly simple form with a known easily-computed inverse leading, effectively, to an explicit update. Such a property is independent of the particular form of the potential energy $V({\bf q})$, provided the non-negativity constraint \eqref{Vpos} is satisfied. 

\subsection{Potential Energy Splitting}
\label{splitting_sec}
In some cases, the expression for the potential energy naturally takes the form
\begin{equation}
\label{split_def}
    V({\bf q}) = \frac{1}{2}{\bf q}^{\intercal}{\bf K}{\bf q}+V'({\bf q})\, ,
\end{equation}
where ${\bf K}$ is a symmetric positive semi-definite $N\times N$ matrix, and where $V'\geq 0$. The first term could represent the stored energy of the system due to linear mechanisms, and $V'$ that due to additional nonlinear effects.  Now, using
\begin{equation}
\label{psidef2}
    V' = \tfrac{1}{2}\psi^2\, ,
\end{equation}
one arrives at a system of equations generalizing \eqref{Hamdyn_g}:
\begin{equation}
\label{Hamdyn_gsplit}
   \dot{{\bf q}} = {\bf M}^{-1}{\bf p}\qquad\qquad  \dot{{\bf p}} = -{\bf K}{\bf q}-\psi{\bf g}\qquad\qquad \dot{\psi} = {\bf g}^{\intercal}\dot{{\bf q}}\, .
\end{equation}
Though equivalent to \eqref{Hamdyn_g}, in a numerical setting, such a splitting allows for a larger family of numerical designs, treating the linear and nonlinear parts of the problem separately. See Section \ref{split_sec}. 

\subsection{Potential Energy Gauge and Regularisation}
\label{shift_sec}
A well-known technique in IEQ/SAV approaches is the introduction of a constant shift of the global energy (see e.g. \cite{SHEN2018407,Liu_SIAM_2020})---Hamilton's equations \eqref{Hamdyn} are unchanged (i.e. gauge invariant), but some numerical schemes, such as IEQ/SAV-based approaches are affected \cite{Lin_JCP_2019}. This amounts to the replacement
\begin{equation}\label{eq:shift}
    V\rightarrow V+\epsilon
\end{equation}
for a suitably chosen shift constant $\epsilon\geq 0$. Under the non-negativity constraint \eqref{Vpos}, $V+\epsilon$ is thus bounded away from zero, regularizing the calculation of ${\bf g}$ as defined in \eqref{gdef}. The regularization approach applies equally to the case of a split potential, under the replacement $V'\rightarrow V'+\epsilon$. 

\section{Numerical Methods}
\label{numsec}
In this section, we assume time discretization using a constant time step $k$, such that solutions are computed at times $t^{n} = nk$, $n=0,1,\hdots$. A discrete-time vector ${\bf u}^{n}$ represents an approximation to a continuous time vector ${\bf u}(t)$ at times $t=t^{n}$. 

\subsection{St{\"o}rmer-Verlet Integration}
\label{Stormersec}

The most basic approach to the integration of the Hamiltonian system as defined in Section \ref{Hamsec} is through direct approximation of the second order form in \eqref{Ham_second}, using centered differences, as:
\begin{equation}
\label{Stormerdef}
    {\bf q}^{n+1} = 2{\bf q}^n-{\bf q}^{n-1}-k^2{\bf M}^{-1}\nabla_{{\bf q}}V|_{{\bf q} = {\bf q}^{n}}\, .
\end{equation}
This discretisation is referred to as St{\"o}rmer-Verlet, and was known to Newton---see \cite{Hairer03} and the references therein. It is a fully explicit two-step scheme---in order to advance the solution to time step $n+1$, a direct evaluation of the gradient of the potential energy $V$ at time step $n$ is required. The values of the state ${\bf q}^{n}$ at $n=0$ and $n=1$ may be initialised, using the initial conditions ${\bf q}^{(0)}$ and ${\bf p}^{(0)}$, as
\begin{equation}
\label{svinit}
    {\bf q}^{0} = {\bf q}^{(0)}\qquad {\bf q}^{1} = {\bf q}^{(0)}+k{\bf M}^{-1}{\bf p}^{(0)}
\end{equation}
to second order in $k$, using a Taylor series approximation. 

St{\"o}rmer-Verlet is second-order accurate in the time step $k$, time-reversible, and symplectic, but not energy-conserving except in an approximate sense \cite{Hairer03}. It can exhibit numerical instability depending on both the choice of time step and the size of the initial conditions. Such instabilities will be mentioned in Section \ref{examplesec}.

\subsection{Explicit Energy-conserving Method (Marazzato et al.)}
\label{Marazzato_sec}
Consider the following scheme approximating \eqref{Hamdyn}, proposed by Marazzato et al. in \cite{Marazzato}:
\begin{equation}
\label{Marazzato}
   {\bf q}^{n+1} = {\bf q}^{n}+k{\bf M}^{-1}{\bf p}^{n+\tfrac{1}{2}}\qquad\qquad {\bf p}^{n+\tfrac{3}{2}} = {\bf p}^{n-\tfrac{1}{2}}-2\int_{nk}^{(n+1)k}\nabla_{{\bf q}}V(\tilde{{\bf q}}^{n}(t)) dt\, .
\end{equation}
This is a time interleaved scheme, with ${\bf p}^{n}$ and ${\bf q}^{n+\tfrac{1}{2}}$ defined at alternating multiples of $k/2$. The first of \eqref{Marazzato} is a standard interleaved approximation to the first of \eqref{Hamdyn}. The second of \eqref{Marazzato} relies on a continuous integration of $\nabla_{{\bf q}}V$ over the time interval $t\in [nk, (n+1)k]$, and over the known free-flight trajectory $\tilde{{\bf q}}^{n}(t)$ defined by
\begin{equation}
\label{free_flight}
\tilde{{\bf q}}^{n}(t) = {\bf q}^{n}+(t-nk){\bf M}^{-1}{\bf p}^{n+\tfrac{1}{2}}\, .
\end{equation}
Other conservative schemes proposed also rely on such a continuous integration \cite{Quispel}. Except for particular functional forms of the potential energy $V({\bf q})$, this integral is not available in closed form and must be approximated. However, once the integral is evaluated, the scheme \eqref{Marazzato} is fully explicit. 

The scheme \eqref{Marazzato} has an associated numerical energy that is exactly conserved \cite{Marazzato}:
\begin{equation}
    \label{Marazzato_egy}
    {\bf H}^{n} = \tfrac{1}{2}{\bf p}^{n+\tfrac{1}{2}}{\bf M}^{-1}{\bf p}^{n-\tfrac{1}{2}}+V({\bf q}^{n}) = {\rm constant}\, .
\end{equation}

There are two important points to mention here:
\begin{itemize}
    \item Scheme \eqref{Marazzato} is exactly conservative, but depends on a continuous integration over a free-flight trajectory as given in \eqref{free_flight} in order to achieve this property. Thus the exact conservation property is approached in the limit of increasing accuracy in the approximation of the continuous integration. It is also important to point out here that a fine-grained approximation of the continuous integration will require multiple evaluations of the gradient $\nabla_{{\bf q}}V$ which, depending on the functional form of $V$, represents an additional cost that grows with the desired accuracy of the approximation. 
    \item When $V\geq 0$, the numerical energy \eqref{Marazzato_egy} does not inherit the non-negativity property of the model system, and thus cannot be used in order to bound solution growth, as pointed out in Remark 8 of \cite{Marazzato}. 
\end{itemize}

Under the additional non-negativity constraint \eqref{Vpos} on $V$, it is possible to demonstrate a method that conserves a pseuodenergy that does not depend on a continuous integration, and that furthermore inherits the non-negativity property of the Hamiltonian of the model system, allowing for useful bounds on solution size. Furthermore, through the exploitation of matrix structure, it will be shown that such a method is fully explicit. Additionally, only one function evaluation $\nabla_{{\bf q}}V$ is required per time step. 

\subsection{Explicit Exactly Energy-conserving Method}

Returning now to the form \eqref{Hamdyn_g} of Hamilton's equations obtained under quadratisation of the potential energy, and the introduction of the new variable $\psi$ as in \eqref{psidef}, consider the following scheme, written in terms of the interleaved time series ${\bf p}^{n+\tfrac{1}{2}}$ and ${\bf q}^{n}$, and the scalar time series $\psi^{n+\tfrac{1}{2}}$:
\begin{subequations}
\label{fd1}
\begin{eqnarray}
    \label{fd1_1}
    {\bf q}^{n+1} &=&{\bf q}^{n}+k{\bf M}^{-1}{\bf p}^{n+\tfrac{1}{2}}\\
    \label{fd1_2}
    {\bf p}^{n+\tfrac{1}{2}} &=& {\bf p}^{n-\tfrac{1}{2}}-\tfrac{k}{2}{\bf g}^{n}\left(\psi^{n+\tfrac{1}{2}}+\psi^{n-\tfrac{1}{2}}\right)\\
    \label{fd1_3}
    \psi^{n+\tfrac{1}{2}} &=& \psi^{n-\tfrac{1}{2}} + \tfrac{1}{2}\left({\bf g}^{n}\right)^{\intercal}\left({\bf q}^{n+1}-{\bf q}^{n-1}\right)\, .
\end{eqnarray}
\end{subequations}
Here, ${\bf g}^{n}$ is defined as
\begin{equation}
\label{gnormal}
 {\bf g}^{n} =  \nabla_{{\bf q}}\psi|_{{\bf q} = {\bf q}^{n}} = \frac{1}{\sqrt{2V({\bf q}^{n})}}\nabla_{{\bf q}}V|_{{\bf q} = {\bf q}^{n}}\, .
\end{equation}

The system \eqref{fd1} may be manipulated into an explicit update form in the following way. Beginning from \eqref{fd1_1}, one may write:
\begin{eqnarray}
{\bf q}^{n+1} &\stackrel{\eqref{fd1_1}}{=}& 2{\bf q}^{n}-{\bf q}^{n-1}+k{\bf M}^{-1}\left({\bf p}^{n+\tfrac{1}{2}}-{\bf p}^{n-\tfrac{1}{2}}\right)\\\notag
 &\stackrel{\eqref{fd1_2}}{=}& 2{\bf q}^{n}-{\bf q}^{n-1}-\tfrac{k^2}{2}{\bf M}^{-1}{\bf g}^{n}\left(\psi^{n+\tfrac{1}{2}}+\psi^{n-\tfrac{1}{2}}\right)\\\notag
   &\stackrel{\eqref{fd1_3}}{=}& 2{\bf q}^{n}-{\bf q}^{n-1}-k^2{\bf M}^{-1}{\bf g}^{n}\psi^{n-\tfrac{1}{2}}-\tfrac{k^2}{4}{\bf M}^{-1}{\bf g}^{n}({\bf g}^{n})^{\intercal}\left({\bf q}^{n+1}-{\bf q}^{n-1}\right)\, .\notag
\end{eqnarray}
Finally, the update has the form
\begin{equation}
\label{linsys}
    {\bf A}^{n}{\bf q}^{n+1} = {\bf b}^{n}\, ,
\end{equation}
where
\begin{equation}
\label{Abdef}
    {\bf A}^{n} = {\bf I}+{\bm \alpha}^{n}({\bm \beta}^{n})^{\intercal}\qquad\qquad {\bf b}^{n} = 2{\bf q}^{n}-2k{\bm \alpha}^{n}\psi^{n-\tfrac{1}{2}}-\left( {\bf I}-{\bm \alpha}^{n}({\bm \beta}^{n})^{\intercal}\right){\bf q}^{n-1}\, .
\end{equation}
Here, ${\bf I}$ is the $N\times N$ identity matrix, and the vectors ${\bm\alpha}^{n}$ and ${\bm\beta}^{n}$ are defined in terms of ${\bf g}^{n}$ by ${\bm\alpha}^{n}= \tfrac{k}{2}{\bf M}^{-1}{\bf g}^{n}$ and  ${\bm\beta}^{n}= \tfrac{k}{2}{\bf g}^{n}$. Thus, given ${\bf q}^{n}$, ${\bf q}^{n-1}$ and $\psi^{n-\tfrac{1}{2}}$, both ${\bf A}^{n}$ and ${\bf b}^{n}$ may be explicitly constructed.  Notice that the matrix ${\bf A}^{n}$ is positive definite by construction, due to the positive definiteness of ${\bf M}$, and thus the update \eqref{linsys} always has a unique solution. Once the update in \eqref{linsys} has been performed, $\psi^{n+\tfrac{1}{2}}$ may be computed from \eqref{fd1_3} explicitly, using ${\bf q}^{n-1}$, ${\bf q}^{n+1}$ and $\psi^{n-\tfrac{1}{2}}$. Scheme \eqref{linsys} requires initial values for ${\bf q}^{0}$ and ${\bf q}^{1}$, which may be set in the same way as for St{\"o}rmer-Verlet, as in \eqref{svinit}, and also $\psi^{\tfrac{1}{2}}$ which may be set to second order in $k$ in terms of the initial displacement ${\bf q}^{(0)}$ and momentum ${\bf p}^{(0)}$ as
\begin{eqnarray}
    \psi^{\tfrac{1}{2}} &=& \sqrt{2V({\bf q}^{(0)})}+\frac{k}{2\sqrt{2V({\bf q}^{(0)})}}\left(\nabla_{{\bf q}}V|_{{\bf q} = {\bf q}^{(0)}}\right)^{\intercal}{\bf M}^{-1}{\bf p}^{(0)}\, .
\end{eqnarray}

Though the solution of \eqref{linsys} apparently requires an $N\times N$ linear system solution, in fact, the inverse is available in closed form, and allows for an explicit solution in $O(N)$ operations. ${\bf A}^{n}$ is a rank-1 perturbation of the identity, and thus Sherman Morrison inversion \cite{Sherman50} yields:
\begin{equation}
    \left({\bf A}^{n}\right)^{-1} = {\bf I}-\frac{{\bm \alpha}^{n}({\bm \beta}^{n})^{\intercal}}{1+({\bm \beta}^{n})^{\intercal}{\bm\alpha}^{n}}\, .
\end{equation}
Computational cost is thus on par with the other methods presented here. Notice in particular, though, that the solution of a linear system involving ${\bf M}^{-1}$ is common to all the methods in this section; the computational cost of performing this operation, which may indeed be the bottleneck, is not considered here. In many cases, though, ${\bf M}$ is diagonal, or even a simple multiple of the identity, and thus its inversion is trivial. 

Exact conservation of a numerical energy follows directly from scheme \eqref{fd1}. Left-multiplying \eqref{fd1_2} by $\tfrac{1}{2}\left({\bf p}^{n+\tfrac{1}{2}}+{\bf p}^{n-\tfrac{1}{2}}\right)^{\intercal}{\bf M}^{-1}$ gives:
\begin{eqnarray}
\tfrac{1}{2}\left({\bf p}^{n+\tfrac{1}{2}}+{\bf p}^{n-\tfrac{1}{2}}\right)^{\intercal}{\bf M}^{-1}\left({\bf p}^{n+\tfrac{1}{2}}-{\bf p}^{n-\tfrac{1}{2}}\right)+\tfrac{k}{4}\left({\bf p}^{n+\tfrac{1}{2}}+{\bf p}^{n-\tfrac{1}{2}}\right)^{\intercal}{\bf M}^{-1}{\bf g}^{n}\left(\psi^{n+\tfrac{1}{2}}+\psi^{n-\tfrac{1}{2}}\right) \stackrel{\eqref{fd1_2}}{=}& 0\\\notag
\tfrac{1}{2}\left({\bf p}^{n+\tfrac{1}{2}}+{\bf p}^{n-\tfrac{1}{2}}\right)^{\intercal}{\bf M}^{-1}\left({\bf p}^{n+\tfrac{1}{2}}-{\bf p}^{n-\tfrac{1}{2}}\right)+\tfrac{1}{4}\left({\bf q}^{n+1}-{\bf q}^{n-1}\right)^{\intercal}{\bf g}^{n}\left(\psi^{n+\tfrac{1}{2}}+\psi^{n-\tfrac{1}{2}}\right) \stackrel{\eqref{fd1_1}}{=}& 0\\\notag
\tfrac{1}{2}\left({\bf p}^{n+\tfrac{1}{2}}+{\bf p}^{n-\tfrac{1}{2}}\right)^{\intercal}{\bf M}^{-1}\left({\bf p}^{n+\tfrac{1}{2}}-{\bf p}^{n-\tfrac{1}{2}}\right)+\tfrac{1}{2}\left(\psi^{n+\tfrac{1}{2}}-\psi^{n-\tfrac{1}{2}}\right)\left(\psi^{n+\tfrac{1}{2}}+\psi^{n-\tfrac{1}{2}}\right) \stackrel{\eqref{fd1_3}}{=}& 0 \, .\notag
\end{eqnarray}
This may be identified as
\begin{equation}
    H^{n+\tfrac{1}{2}} =  H^{n-\tfrac{1}{2}} = {\rm constant}\, ,
\end{equation}
where the numerical energy $H^{n+\tfrac{1}{2}}$ is
\begin{equation}
\label{explicit_energy}
    H^{n+\tfrac{1}{2}} = \tfrac{1}{2}\left({\bf p}^{n+\tfrac{1}{2}}\right)^{\intercal}{\bf M}^{-1}{\bf p}^{n+\tfrac{1}{2}}+\tfrac{1}{2}\left(\psi^{n+\tfrac{1}{2}}\right)^2 \geq 0\, .
\end{equation}
It is worth comparing this expression, for which non-negativity is ensured, with the expression \eqref{Marazzato_egy} for the scheme in Section \ref{Marazzato_sec}. The momentum ${\bf p}^{n+\tfrac{1}{2}}$ may be bounded in terms of the energy $H$, for all time steps $n$, as
\begin{equation}
    \|{\bf p}^{n+\tfrac{1}{2}}\|\leq \sqrt{2\lambda_{{\rm max}}({\bf M})H}\, ,
\end{equation}
which is identical to the bound \eqref{pbound} for the model system. The scheme \eqref{fd1} is thus unconditionally numerically stable. 

\subsection{Split Potential Form}
\label{split_sec}

Consider now the split form of the potential energy $V$, as described in Section \ref{splitting_sec}, where a positive semi-definite quadratic form has been separated from $V$ to leave a residual energy contribution $V'\geq 0$, from which an auxiliary variable $\psi$ may be defined as in \eqref{psidef2}. Consider the following scheme, now modified with respect to \eqref{fd1} through the addition of a linear term in \eqref{fd2_2} below:
\begin{subequations}
\label{fd2}
\begin{eqnarray}
    \label{fd2_1}
    {\bf q}^{n+1} &=&{\bf q}^{n}+k{\bf M}^{-1}{\bf p}^{n+\tfrac{1}{2}}\\
    \label{fd2_2}
    {\bf p}^{n+\tfrac{1}{2}} &=& {\bf p}^{n-\tfrac{1}{2}}-k{\bf K}{\bf q}^{n}-\tfrac{k}{2}{\bf g}^{n}\left(\psi^{n+\tfrac{1}{2}}+\psi^{n-\tfrac{1}{2}}\right)\\
    \label{fd2_3}
    \psi^{n+\tfrac{1}{2}} &=& \psi^{n-\tfrac{1}{2}} + \tfrac{1}{2}\left({\bf g}^{n}\right)^{\intercal}\left({\bf q}^{n+1}-{\bf q}^{n-1}\right)\, .
\end{eqnarray}
\end{subequations}
where now,
\begin{equation}
\label{gsplit}
 {\bf g}^{n} =  \nabla_{{\bf q}}\psi|_{{\bf q} = {\bf q}^{n}} = \frac{1}{\sqrt{2V'({\bf q}^{n})}}\nabla_{{\bf q}}V'|_{{\bf q} = {\bf q}^{n}}\, .
\end{equation}
An explicit update follows as in the case of the non-split form in \eqref{linsys}, with ${\bf A}^{n}$ as given in \eqref{Abdef}, but with ${\bf b}^{n}$ now defined as
\begin{equation}
{\bf b}^{n} = \left(2{\bf I}-k^2 {\bf M}^{-1}{\bf K}\right){\bf q}^{n}-2k{\bm \alpha}^{n}\psi^{n-\tfrac{1}{2}}-\left( {\bf I}-{\bm \alpha}^{n}({\bm \beta}^{n})^{\intercal}\right){\bf q}^{n-1}\, .
\end{equation}

An expression for a conserved numerical energy follows as before, now taking the form:
\begin{equation}
\label{Egy_split}
    H^{n+\tfrac{1}{2}} = \tfrac{1}{2}\left({\bf p}^{n+\tfrac{1}{2}}\right)^{\intercal}{\bf M}^{-1}{\bf p}^{n+\tfrac{1}{2}} +\tfrac{1}{2}({\bf q}^{n+1})^{\intercal}{\bf K}{\bf q}^{n}+ \tfrac{1}{2}\left(\psi^{n+\tfrac{1}{2}}\right)^2 = {\rm constant}\, .
\end{equation}
Because the second term is of indefinite sign, the numerical energy is not necessarily non-negative. Because it is a quadratic form, however, it may be bounded, using \eqref{fd2_1}, as:
\begin{equation}
    \tfrac{1}{2}({\bf q}^{n+1})^{\intercal}{\bf K}{\bf q}^{n}\geq -\tfrac{1}{8}({\bf q}^{n+1}-{\bf q}^{n})^{\intercal}{\bf K}({\bf q}^{n+1}-{\bf q}^{n}) = -\tfrac{k^2}{8}\left({\bf p}^{n+\tfrac{1}{2}}\right)^{\intercal}{\bf M}^{-1}{\bf K}{\bf M}^{-1}{\bf p}^{n+\tfrac{1}{2}}\, .
\end{equation}
It then follows that
\begin{equation}
    H^{n+\tfrac{1}{2}} \geq \tfrac{1}{2}\left({\bf p}^{n+\tfrac{1}{2}}\right)^{\intercal}\left({\bf M}^{-1}-\frac{k^2}{4}{\bf M}^{-1}{\bf K}{\bf M}^{-1}\right){\bf p}^{n+\tfrac{1}{2}}\, .
\end{equation}
Given that ${\bf M}>0$, ${\bf M}^{-1}$ may be factored into unique upper and lower triangular factors as ${\bf M}^{-1} = {\bf M}^{-\tfrac{T}{2}}{\bf M}^{-\tfrac{1}{2}}$, a condition for non-negativity of $H$ is
\begin{equation}
\label{stabcond}
    k\leq \frac{2}{\lambda_{{\rm max}}\left({\bf M}^{-\tfrac{1}{2}}{\bf K}{\bf M}^{-\tfrac{T}{2}}\right)}\, .
\end{equation}
The scheme is now conditionally stable, with the stability condition \eqref{stabcond} corresponding to that for a linear system when $V' = 0$. 

\subsection{Remarks}
\label{remarkssec}
\begin{itemize}
    \item {\bf Linear Conditions:} One of the interesting features of the non-split form \eqref{fd1} is that, even if the potential $V$ is quadratic in ${\bf q}$, implying a linear system, the scheme is not linear. If, however, the split form \eqref{fd2} is used, then under linear conditions, $V'=0$, and scheme \eqref{fd2} reduces exactly to St{\"o}rmer-Verlet \eqref{Stormerdef}. Note that under linear conditions, St{\"o}rmer-Verlet does indeed possess an exactly conserved numerical energy, as given in \eqref{Egy_split} with $\psi^{n+\tfrac{1}{2}} = 0$. 
    \item {\bf Generalized Splitting:} The non-split form \eqref{fd1} and the split form \eqref{fd2} may be combined to yield a large family of conservative schemes in an obvious way. If 
    \begin{equation}
        V({\bf q}) = \tfrac{1}{2}{\bf q}^{\intercal}{\bf K}{\bf q}+V_{{\rm nonlinear}}({\bf q})\, ,
    \end{equation}
    where $V_{{\rm nonlinear}}\geq 0$ and consists of higher order terms in ${\bf q}$, then any splitting of the form
    \begin{equation}
        V({\bf q}) = \tfrac{1}{2}{\bf q}^{\intercal}{\bf K}_{0}{\bf q}+V'({\bf q})\qquad V'({\bf q}) = V_{{\rm nonlinear}}+\tfrac{1}{2}{\bf q}^{\intercal}{\bf K}'{\bf q}\, ,
    \end{equation}
    with ${\bf K} = {\bf K}_{0}+{\bf K}'$, ${\bf K}_{0}\geq 0$, ${\bf K}'\geq 0$ will yield a conservative form. The resulting numerical scheme, with $V' = \tfrac{1}{2}\psi^2$ will inherit conservation of energy, which will be non-negative under a condition analogous to \eqref{stabcond}, but depending on ${\bf K}_{0}$. 
    \item {\bf Generalized Update for ${\bf q}^{n}$:} Consider again the non-split scheme \eqref{fd1}, which may be rewritten as
    \begin{subequations}
\label{fd4}
\begin{eqnarray}
    \label{fd4_1}
    {\bf q}^{n+1} &=&{\bf q}^{n}+k{\bf M}^{-1}{\bf p}^{n+\tfrac{1}{2}}\\
    \label{fd4_2}
    {\bf p}^{n+\tfrac{1}{2}} &=& {\bf p}^{n-\tfrac{1}{2}}-\tfrac{k}{2}{\bf g}^{n}\left(\psi^{n+\tfrac{1}{2}}+\psi^{n-\tfrac{1}{2}}\right)\\
    \label{fd4_3}
    \psi^{n+\tfrac{1}{2}} &=& \psi^{n-\tfrac{1}{2}} + \tfrac{k}{2}\left({\bf g}^{n}\right)^{\intercal}{\bf M}^{-1}\left({\bf p}^{n+\tfrac{1}{2}}+{\bf p}^{n-\tfrac{1}{2}}\right)\, .
\end{eqnarray}
\end{subequations}
In this form, exact energy conservation follows from \eqref{fd4_2} and \eqref{fd4_3} only---it is independent of the values of ${\bf g}^{n}$, which are derived solely from ${\bf q}^{n}$. Further opportunities for generalization are thus available---\eqref{fd4_1} could be replaced by any consistent update for ${\bf q}^{n+1}$, and the exact numerical energy conservation property remains undisturbed. 

\item {\bf Variable Time Steps:} Though the case of variable time steps will not be discussed here in detail (and was indeed investigated in \cite{Marazzato}), an exact energy-conserving scheme follows immediately from the form given in \eqref{fd4} above. Consider time instants $t^{n}$, and $t^{n+1/2}$, for integer $n$, where $t^{n}<t^{n+1/2}<t^{n+1}$. From these, one may define two sequences of time steps: $k^{n} = t^{n+1/2}-t^{n-1/2}$ and $k^{n+1/2} = t^{n+1}-t^{n}$. Supposing that ${\bf q}^{n}$ is a time series defined for $t = t^{n}$, and similarly ${\bf p}^{n+\tfrac{1}{2}}$ and ${\bf \psi}^{n+\tfrac{1}{2}}$ are defined for $t=t^{n+\tfrac{1}{2}}$, then a scheme follows as:
 \begin{subequations}
\label{fd5}
\begin{eqnarray}
    \label{fd5_1}
    {\bf q}^{n+1} &=&{\bf q}^{n}+k^{n+\tfrac{1}{2}}{\bf M}^{-1}{\bf p}^{n+\tfrac{1}{2}}\\
    \label{fd5_2}
    {\bf p}^{n+\tfrac{1}{2}} &=& {\bf p}^{n-\tfrac{1}{2}}-\tfrac{k^{n}}{2}{\bf g}^{n}\left(\psi^{n+\tfrac{1}{2}}+\psi^{n-\tfrac{1}{2}}\right)\\
    \label{fd5_3}
    \psi^{n+\tfrac{1}{2}} &=& \psi^{n-\tfrac{1}{2}} + \tfrac{k^{n}}{2}\left({\bf g}^{n}\right)^{\intercal}{\bf M}^{-1}\left({\bf p}^{n+\tfrac{1}{2}}+{\bf p}^{n-\tfrac{1}{2}}\right)\, .
\end{eqnarray}
\end{subequations}
It is direct to show that, regardless of the choices of $t^{n}$ and $t^{n+\tfrac{1}{2}}$, the scheme \eqref{fd5} conserves the energy \eqref{explicit_energy} exactly. Thus such a generalisation is unconditionally stable, under the same reasoning as in the case of constant time steps. 

\item {\bf Regularisation:} Regularisation through a potential energy shift, as described in Section \ref{shift_sec}, impacts on the calculation of ${\bf g}$ in \eqref{gnormal} and \eqref{gsplit}, through a replacement $V\rightarrow V+\epsilon$ or $V'\rightarrow V'+\epsilon$, respectively. 
\end{itemize}

\section{Examples}
\label{examplesec}
In this section, various Hamiltonian systems are simulated using schemes \eqref{fd1} and \eqref{fd2}, beginning with the Fermi-Pasta-Ulam ODE problem in Section \ref{fpusec}, and then progressing to more complex ODE systems derived as semi-discretisations to PDE systems. These include the coupled transverse-longitudinal vibration of a string at high amplitudes, in Section \ref{stringsec}, and the high amplitude vibration of a thin plate in Section \ref{platesec}. 

In all cases, a useful measure of the exact energy conservation property is the relative energy deviation error, 
\begin{equation}
\label{relerrdef}
   {\rm  Relative\,\, error} = \frac{H^{n+\tfrac{1}{2}}-H^{\tfrac{1}{2}}}{H^{\tfrac{1}{2}}}\, .
\end{equation}
In double precision floating point arithmetic, it is normally on the order of machine precision, or approximately 10$^{-16}$. Depending on the state size of the system in question, however, larger deviations are possible. 

\subsection{Fermi-Pasta-Ulam Problem}
\label{fpusec}
As a simple example, consider the classic system of a linear arrangement masses connected by linear and nonlinear springs, as proposed originally by Fermi, Pasta and Ulam \cite{Fermi}, and later adapted as a test problem by various authors \cite{Hairer, Marazzato}, the form of which is followed here. 

Consider a system of $N=2M$ masses in a linear arrangement, which longitudinal displacements $q_{i}$ and momenta $p_{i}$, $i=1,\hdots,2M$. The system Hamiltonian is defined by
\begin{equation}
\label{fpudef}
    H = \frac{1}{2}\sum_{i=1}^{2M}p_{i}^{2}+V\qquad {\rm where}\qquad V = \frac{\omega^2}{4}\sum_{i=1}^{M}\left(q_{2i}-q_{2i-1}\right)^2+\sum_{i=0}^{M}\left(q_{2i+1}-q_{2i}\right)^4\, ,
\end{equation}
where, in the expression above, $q_{0} = q_{2M+1} = 0$. Thus each mass is connected, in an alternating arrangement, to a linear spring, and a cubic nonlinear spring. 

$V$ is clearly non-negative here, and thus, after consolidation of the displacements $q_{i}$ and $p_{i}$ into vectors ${\bf p}$ and ${\bf q}$ of size $N\times 1$, the scheme as presented in \eqref{fd1} follows immediately, with ${\bf M} = {\bf I}_{N}$, the $N\times N$ identity matrix, and with $V$ as given in \eqref{fpudef} above. This scheme is explicit, exactly conservative, and unconditionally stable. 

A split form \eqref{split_def} of the Hamiltonian also follows, using
\begin{equation}
    {\bf K} = \frac{\omega^2}{2} {\bf I}_{M}\otimes\begin{bmatrix}
    \,\,\,1 & \!\!\!-1\\
    -1 & 1 \\
    \end{bmatrix}\qquad\qquad V' = \sum_{i=0}^{M}\left(q_{2i+1}-q_{2i}\right)^4 \geq 0\, ,
\end{equation}
where here, $\otimes$ indicates a Kronecker product. An exactly energy conserving method follows, as in \eqref{fd2}, and is stable under the condition \eqref{stabcond}, which reduces in this case to \begin{equation}
    k\leq \frac{2}{\omega}\, .
\end{equation}

\subsubsection{Numerical Results}

\begin{figure}[hbt!]
\centering
\includegraphics[width = 0.9\linewidth]{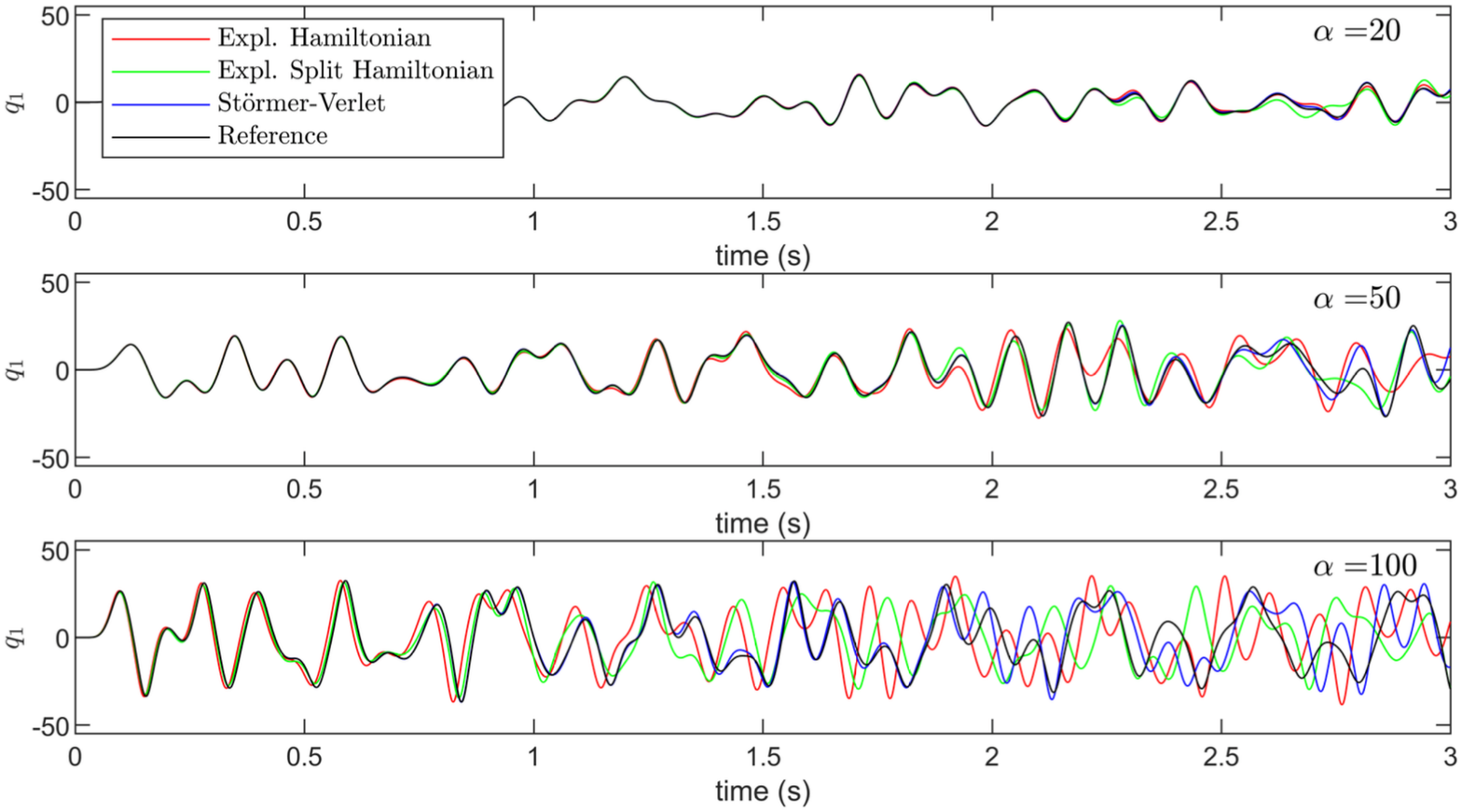}
\caption{Time histories of $q_{1}$, under the initial condition $q_{4} = \alpha$, and under different values of $\alpha$ as indicated. Results for schemes \eqref{fd1}, \eqref{fd2} and \eqref{Stormerdef} are plotted against a reference solution. }
\label{fpuevfig}
\end{figure}

\begin{figure}[hbt!]
\centering
\includegraphics[width = 0.9\linewidth]{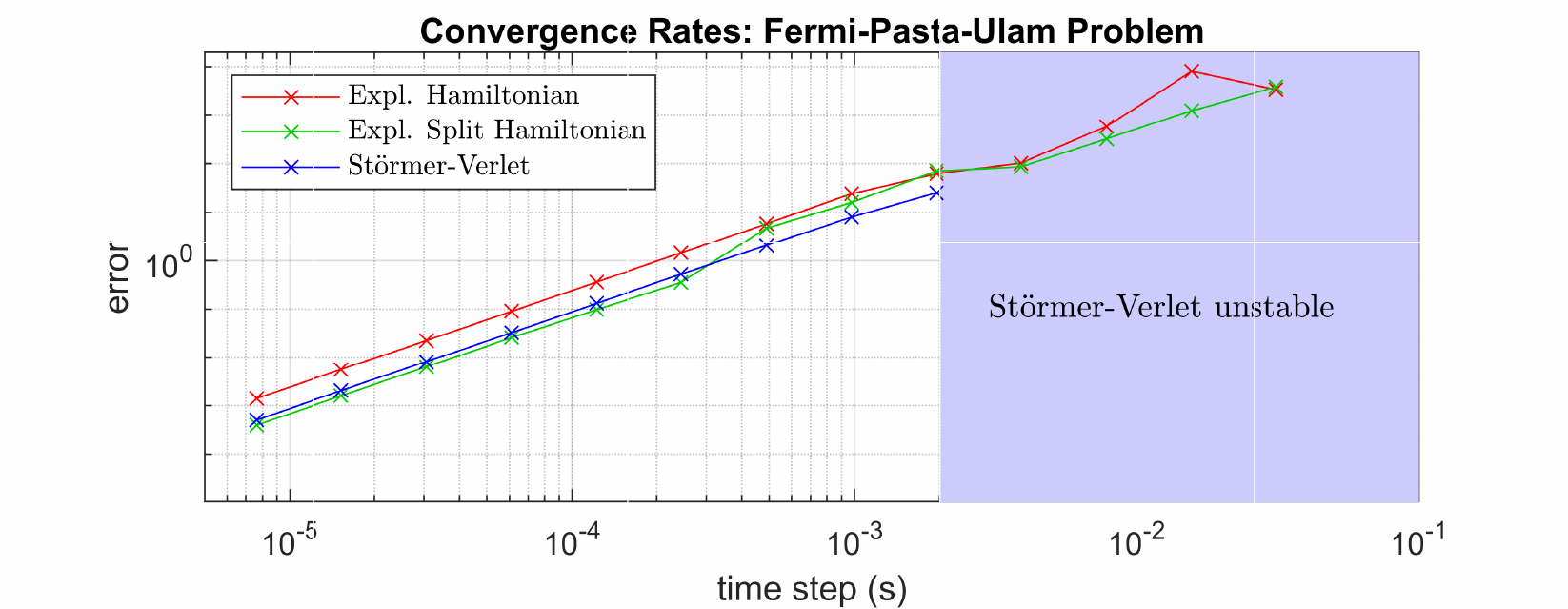}
\caption{Error, as defined in \eqref{errordef}, as a function of time step $k$ for the explicit Hamiltonian scheme \eqref{fd1}, the split potential form \eqref{fd2}, and St{\"o}rmer-Verlet \eqref{Stormerdef}, for the Fermi-Pasta-Ulam problem, where ${\bf q}^{(0)} = [0,0,0,100,0,0]$. }
\label{fpuerrorfig}
\end{figure}

We use the same settings as in \cite{Marazzato}, and choose $\omega = 50$, and $M = 3$. Simulations are run here with a time step of $k = 10^{-3}$ s; the reference solution is computed using St{\"o}rmer-Verlet, with a time step of $k=2^{-20}\approxeq 10^{-6}$ s. For initial conditions, we set ${\bf p}^{(0)} = {\bf 0}$, and ${\bf q}^{(0)} = [0,0,0,\alpha,0,0]^{\intercal}$, for different values of $\alpha$. See Figure \ref{fpuevfig}, illustrating c comparison between time histories using the exactly conservative schemes \eqref{fd1} and \eqref{fd2} with St{\"o}rmer-Verlet \eqref{Stormerdef}. In all cases, the onset of errors is slightly faster for \eqref{fd1} than for the other schemes, indicating a higher error (see also Figure \ref{fpuerrorfig}). 

\begin{figure}[ht]
\centering
\includegraphics[width = 0.9\linewidth]{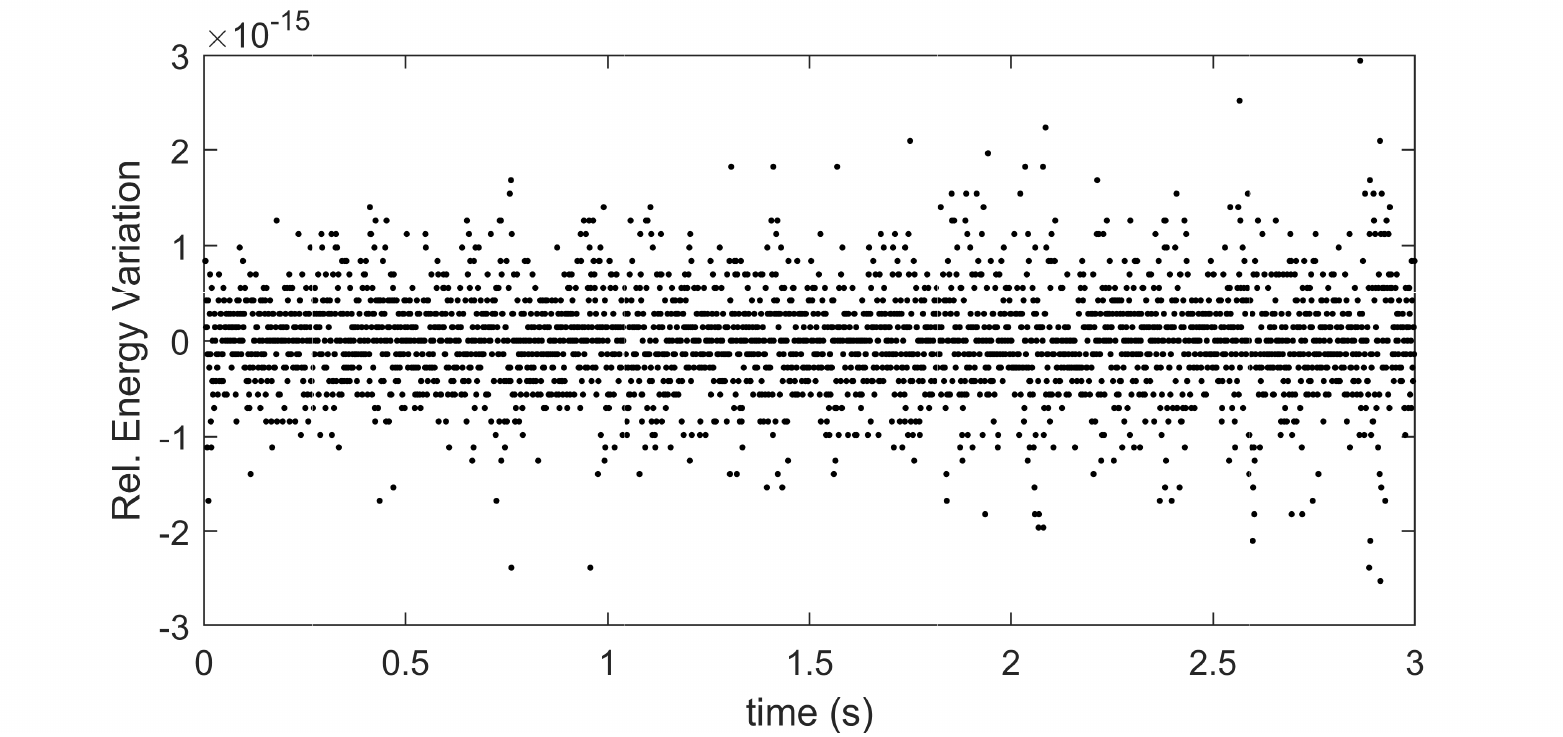}
\caption{Relative variation in numerical energy, as defined by \eqref{relerrdef}, for scheme \eqref{fd1} for the Fermi-Pasta-Ulam system. }
\label{fpuegyfig}
\end{figure}

For the schemes \eqref{fd1} and \eqref{fd2}, numerical energy, as defined by \eqref{explicit_energy} and \eqref{Egy_split} respectively, is conserved to machine precision. See Figure \ref{fpuegyfig}, illustrating the relative deviation in energy, as defined in \eqref{relerrdef}, for scheme \eqref{fd1} for the Fermi-Pasta-Ulam system, under the conditions as described above, and for $\alpha = 100$. It is easily seen that the relative error is on the order of machine roundoff error in double-precision floating point, or approximately $10^{-16}$. Furthermore, it is possible to directly observe the quantisation of the relative energy to the lowest bits in the machine number representation. 

As a basic test of convergence, consider the $L^{2}$ error defined, over a simulation duration $n=0,\hdots,N_{f}$, by
\begin{equation}
\label{errordef}
    {\rm Error} = \sqrt{\sum_{n=0}^{N_{f}}k\|{\bf q}^{n}-{\bf q}_{{\rm ref}}(t=nk)\|^2}\, ,
\end{equation}
where ${\bf q}_{{\rm ref}}$ are computed using St{\"o}rmer-Verlet with a time step of $k=2^{-20}\approxeq 10^{-6}$ s. The error appears in Figure \ref{fpuerrorfig}, for the exact energy conserving method in \eqref{fd1}, the split potential method in \eqref{fd2}, and using St{\"o}rmer-Verlet, over a range of time steps. In this case, the total simulation duration is 1 s, and the initial condition is chosen to be large, with  $\alpha = 100$. No regularisation was employed in this case (see Section \ref{shift_sec}). In general, the errors for St{\"o}rmer-Verlet and the split scheme \eqref{fd2} track each other quite closely, with the non-split method \eqref{fd1} performing somewhat worse. Second order accuracy is easily observed in all cases. Stability for St{\"o}rmer-Verlet is dependent on the size of the initial condition; for $\alpha = 100$, the range of time steps over which St{\"o}rmer-Verlet is unstable is indicated in the figure. 

\subsection{Nonlinear String Vibration}
\label{stringsec}

As a second example, consider the motion of a taut string. For sufficiently large displacements, the linear wave equation must be augmented by appropriate nonlinear terms to account for the amplitude-dependent physical phenomena observed during motion. 
One may obtain a  geometrically exact model by considering large strains, and applying Hooke's law \cite{morse_book}. This model is written compactly as
\begin{equation}\label{eq:NlinStringCnt}
    \rho A \partial_t^2 { u}=\partial_x \left(\frac  {\partial {\mathcal V}}{\partial \zeta}\right), \quad  \rho A \partial_t^2 { v}=\partial_x \left(\frac  {\partial {\mathcal V}}{\partial \eta}\right).
\end{equation}
Here, $u=u(x,t) : \mathcal{D} \times \mathbb R_0^+ \rightarrow \mathbb{R}$ represents the transverse displacement, for a spatial domain $\mathcal{D} =[0,L]$, where $L$ is the length of the unstretched string. Similarly, $v = v(x,t)$ is the longitudinal displacement. In \eqref{eq:NlinStringCnt}, $\partial_{t}$ and $\partial_{x}$ represent partial differentiation with respect to time $t$ and spatial coordinate $x$. Furthermore, we define $\zeta = \partial_x u$, $\eta = \partial_x v$. The function ${\mathcal V} = {\mathcal V}(\zeta, \eta): \mathbb{R}^2 \rightarrow \mathbb{R}_0^+$ is a non-negative potential density, which may be given in two equivalent forms as
\begin{subnumcases}{ {\mathcal V}(\zeta, \eta) = \label{eq:NlinStringVCnt}}
 \frac{EA}{2}\left(\zeta^2 + \eta^2 \right) - (EA-T_0) \left( \sqrt{(1+\eta)^2 + \zeta^2} -1-\eta\right), \label{eq:NlinStringVCnt1} \\
 \frac{T_0}{2}\left(\zeta^2 + \eta^2 \right) + \frac{EA-T_0}{2} \left( \sqrt{(1+\eta)^2 + \zeta^2} -1\right)^2. \label{eq:NlinStringVCnt2}
\end{subnumcases}
Here, the various constants that appear are: $\rho$, the volume density in kg$\cdot$m$^{-3}$; $E$, Young's modulus, in kg$\cdot$ s$^{-2}$ $\cdot$ m$^{-1}$; $A$, the area of the string cross section in m$^{2}$; and $T_0$, the applied tension in kg$\cdot$ m $\cdot$ s$^{-2}$. It is easy to show that \eqref{eq:NlinStringVCnt1} is the same as \eqref{eq:NlinStringVCnt2}, and that ${\mathcal V} \geq 0 \,\, \forall \,\, (\zeta,\, \eta)$. In \eqref{eq:NlinStringVCnt2}, one may easily  split $\mathcal V$ into a quadratic part, plus a non-negative nonlinear term, provided that $EA > T_0$ (a condition that is generally satisfied, for instance by all strings of interest in musical acoustics): ultimately, this allows for numerical solution using a split potential form as described in Section \ref{split_sec}. 
Both forms have been used in previous works: \eqref{eq:NlinStringVCnt1} in e.g. \cite{chabassier_CMAME_2010,Marazzato};  \eqref{eq:NlinStringVCnt2} in \cite{ducceschi_JSV_2021}.
System \eqref{eq:NlinStringCnt} is Hamiltonian, with the total energy defined by
\begin{equation}
    H = \int_{\mathcal D} \frac{\rho A}{2}\, \left( (\partial_t u)^2 + (\partial_t v)^2 \right) + \mathcal{V}\, \mathrm{d}x.
\end{equation}
Energy conservation holds under a suitable set of boundary conditions. Here, conditions of fixed type are considered, such that $u=v=0 $ at $x=\{0,L\}$.

\subsubsection{Semi-discrete Form}

The domain $\mathcal D$ may be divided into segments of length $h$, the grid spacing. Let $M$ be the total number of segments, yielding $M-1$ grid points, not including the end points, where the solution is fixed to zero. The continuous functions $u(x,t)$, $v(x,t)$ may then be approximated by grid functions $u_l(t)$, $v_l(t)$, at the grid point $x_l = lh$, $l=1,...,M-1$. Approximations to $\partial_x$ may be given as difference operators, expressed here in terms of their action on a grid function $u_l(t)$:
\begin{equation}\label{eq:Dxp1D}
    D_{+} u_l = \frac{1}{h}\left( u_{l+1} - u_{l}\right), \,\,\,   D_{-} u_l = \frac{1}{h}\left( u_{l} - u_{l-1}\right).
\end{equation}
From these, one may also define the second difference operator as $D_2 = D_+ D_-$.  Furthermore, let $\zeta_l = D_-u_l$, $\eta_l = D_-v_l$. A semi-discrete form for \eqref{eq:NlinStringCnt} is then obtained as
\begin{equation}\label{eq:StringSemiDisc}
\rho A \ddot u_l = D_+ \left(\frac{\partial {\mathcal V}_l}{\partial \zeta_l} \right), \qquad \rho A \ddot v_l = D_+ \left(\frac{\partial {\mathcal V}_l}{\partial \eta_l} \right),
\end{equation}
where ${\mathcal V}_l \triangleq {\mathcal V}(\zeta_l,\eta_l)$. It is convenient to write system \eqref{eq:StringSemiDisc} compactly using the consolidated state vector ${\bf q} = [{\bf u}^\intercal,{\bf v}^\intercal]^\intercal$. In this form, the difference operators are represented by matrices, such that $D_-$ becomes the $M \times (M-1)$ matrix ${\bf D}_-$ given in terms of its action on e.g. $\bf u$ as
\begin{equation}
  {\bf D}_- {\bf u}= \frac{1}{h}([{\bf u}^\intercal,0] - [0, {\bf u}^\intercal]).
\end{equation}
Then,  define ${\bf D}_+ = -{\bf D}_-^\intercal$, and ${\bf D}_2 = {\bf D}_+ {\bf D}_-$. In vector form, \eqref{eq:StringSemiDisc} becomes
\begin{equation}\label{eq:StringSemiDiscq}
    \rho A \ddot {\bf q} = \frac{1}{h}\begin{bmatrix}{\bf D}_+\nabla_{\boldsymbol \zeta}{V} \\ {\bf D}_+\nabla_{\boldsymbol \eta}{V} \end{bmatrix},
\end{equation}
where ${\boldsymbol \zeta} = {\bf D}_- {\bf u}$, ${\boldsymbol \eta} = {\bf D}_- {\bf v}$. This system conserves the semi-discrete energy of the form \eqref{Hamdef}, with 
\begin{equation}\label{eq:EnStringSemiDisc}
   {\bf M} = \rho A h {\bf I}_{2M-2}, \quad V = h \sum_{l=1}^M {\mathcal V}_l.
\end{equation}
Here ${\bf I}_{2M-2}$ is the $(2M-2) \times (2M-2)$ identity matrix. A proof is obtained immediately by noting that
\begin{equation}
    \dot V = \dot{\boldsymbol \zeta}^\intercal \nabla_{\boldsymbol \zeta}V + \dot {\boldsymbol \eta}^\intercal \nabla_{\boldsymbol \eta}V = -\dot {\bf u}^\intercal{\bf D}_+\nabla_{\boldsymbol \zeta}V - \dot {\bf v}^\intercal{\bf D}_+\nabla_{\boldsymbol \eta}V = -\dot {\bf q}^\intercal \begin{bmatrix}{\bf D}_+\nabla_{\boldsymbol \zeta}{V} \\ {\bf D}_+\nabla_{\boldsymbol \eta}{V} \end{bmatrix}.
\end{equation}
Thus, left-multiplying \eqref{eq:StringSemiDiscq}  by $h \dot {\bf q}^\intercal$ yields \eqref{eq:EnStringSemiDisc}. A split potential form is also available, via \eqref{eq:NlinStringVCnt2}. In this case,
\begin{equation}\label{eq:EnStringSemiDiscSplit}
    {\bf K} = - T_0 h \begin{bmatrix}{\bf D}_2 & {\bf 0} \\ {\bf 0} & {\bf D}_2\end{bmatrix}, \quad V^\prime = \frac{h(EA-T_0)}{2} \sum_{l=1}^M   \left( \sqrt{(1+\eta_l)^2 + \zeta_l^2} -1\right)^2
\end{equation}

\subsubsection{Numerical Methods}
System \eqref{eq:StringSemiDiscq} may be integrated in time in a number of ways. First, introduce the discrete time vector ${\bf q}^n$, approximating ${\bf q}(t)$ at the time $t_n=kn$. Furthermore, let $V^n = V({\boldsymbol \zeta}^n,{\boldsymbol \eta}^n)$. An explicit time stepping scheme is obtained by application of the St{\"o}rmer-Verlet algorithm:
\begin{equation}\label{eq:SVstring}
    {\bf q}^{n+1} = 2 {\bf q}^{n+1} - {\bf q}^{n-1} + \frac{k^2}{\rho A h} \begin{bmatrix}{\bf D}_+\nabla_{{\boldsymbol \zeta}^n}{V}^n \\ {\bf D}_+\nabla_{{\boldsymbol \eta}^n}{V}^n \end{bmatrix}.
\end{equation}
While simple, this scheme does not conserve a positive discrete energy, and instabilities may occur at large displacements.  A stable scheme may be obtained by an energy-conserving discretisation of the gradient \cite{itoh_JCP_1988,ducceschi:ICA:2019}:
\begin{equation}\label{eq:NTstring}
     {\bf q}^{n+1} = 2 {\bf q}^{n+1} - {\bf q}^{n-1} +  \frac{k^2}{\rho A } \begin{bmatrix}{\bf D}_+ \, \mathfrak{g}^n_{{\boldsymbol \zeta}} \\ {\bf D}_+ \,   \mathfrak{g}^n_{{\boldsymbol \eta}}\end{bmatrix}
\end{equation}
where the discrete gradients are defined as
\begin{equation}
    (\mathfrak{g}^n_{{\boldsymbol \zeta}})_l = \frac{\mathcal{V}({ \zeta}^{n+1}_l,{ \eta}^{n}_l)-\mathcal{V}({ \zeta}_l^{n-1},{ \eta}_l^{n})}{{ \zeta}_l^{n+1}-{ \zeta}^{n-1}_l}, \quad    (\mathfrak{g}^n_{{\boldsymbol \eta}})_l = \frac{\mathcal{V}({ \zeta}^{n}_l,{ \eta}^{n+1}_l)-\mathcal{V}({ \zeta}_l^{n},{ \eta}_l^{n-1})}{{ \eta}_l^{n+1}-{ \eta}^{n-1}_l},
\end{equation}
with $l = 1,..,M$. This scheme leads to discrete energy conservation, and unconditional stability, but it is fully implicit, and will generally require the use of iterative root finding routines such as Newton-Raphson \cite{ducceschi:ICA:2019}. 

Finally, scheme \eqref{fd1} results from Hamiltonian \eqref{Hamdef} with definitions \eqref{eq:EnStringSemiDisc}. The split-potential form \eqref{fd2} is also available, via \eqref{eq:EnStringSemiDiscSplit}. In the latter case, a stability condition arises as per \eqref{stabcond}, such that
\begin{equation}
    k \leq \sqrt{\rho A / T_0} \, \, h.
\end{equation}
This is the standard Courant-Friedrichs-Lewy stability condition for the one-dimensional linear wave equation \cite{Courant28}. 
\subsubsection{Numerical Results}

\begin{figure}[hbt!]
    \centering
    \includegraphics[width = 0.9\linewidth]{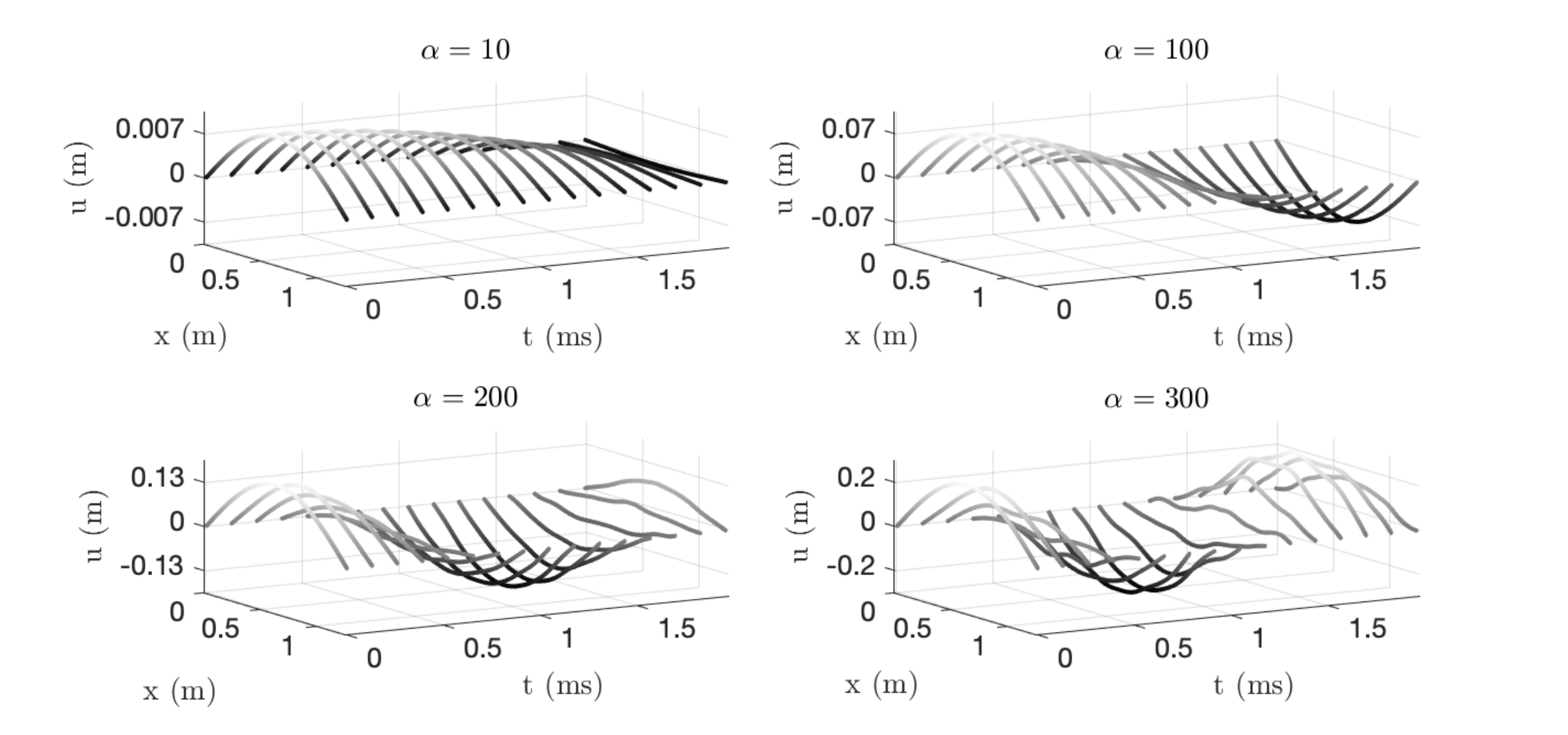}
    \caption{Snapshots of the geometrically exact nonlinear string, computed using scheme \eqref{fd2},  with split potential form as per \eqref{eq:EnStringSemiDiscSplit}. Here, $k = 2.6 \cdot 10^{-7}$, and the grid spacing is chosen as $h = 1.05 \sqrt{E/\rho}\, k$. 
    The string's initial normalised amplitude $\alpha$ is as indicated.}
    \label{fig:snapsC3string}
\end{figure}

As a first example, let the string be initialised in its first linear mode of vibration in the transverse direction, that is
\begin{equation}
    u(x,0) = \alpha \sqrt{A} \sin(\pi x/L), \quad v(x,0) = 0,
\end{equation}
and let the initial velocity be zero for both transverse and longitudinal motion. The amplitude parameter $\alpha$ is nondimensional. The string parameters are taken from \cite{chabassierPianoParams}, for the C3 piano string, and are: $\rho = 7850$ kg m$^{-3}$; $A=8.87\cdot 10^{-7}$ m$^2$; $L = 1.259$ m; $E = 2.02 \cdot 10^{11}$ kg$\cdot$ s$^{-2}$ $\cdot$ m$^{-1}$; $T_0 = 759$ kg$\cdot$ m $\cdot$ s$^{-2}$. Figure \ref{fig:snapsC3string} shows snapshots of the computed solution using scheme \eqref{fd2}, and under various initial amplitudes $\alpha$. In this case, a regularisation of the potential energy $V'$ is employed (see Section \ref{shift_sec}), with a shift of $\epsilon = 10^{8}$.  Typical amplitude-dependent phenomena are visible: the frequency of vibration increases with the initial amplitude, and the initial shape deforms progressively. 

\begin{figure}[hbt!]
    \centering
    \includegraphics[width = 0.9\linewidth]{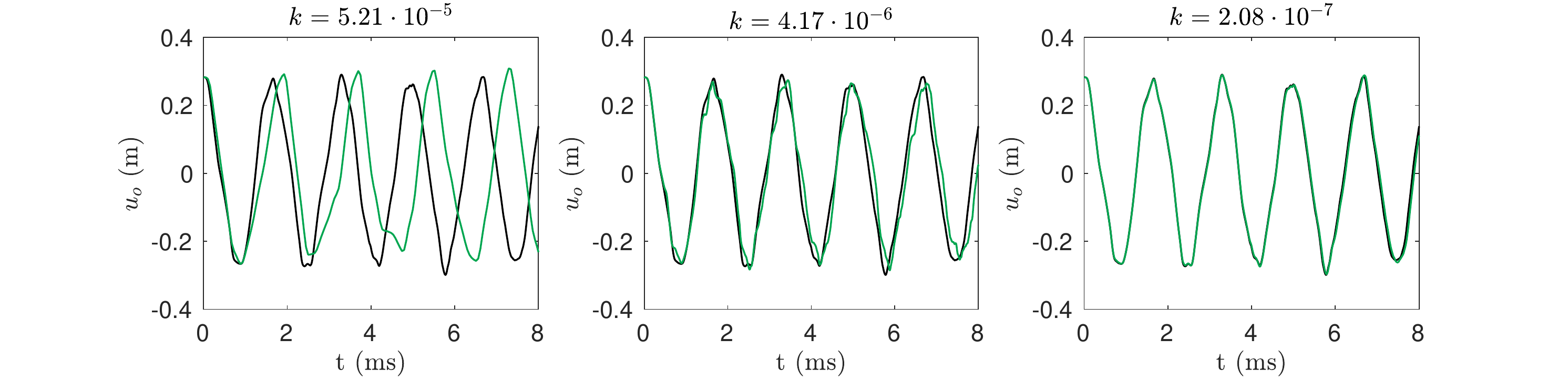}
    \caption{Transverse displacement $u$ (in green) for a string initialised with $\alpha = 300$, at $x = 0.5L$. In black is the reference solution, computed using scheme \eqref{eq:SVstring} with a time step $k=1.04 \cdot 10^{-7}$. In green are the waveforms computed using \eqref{fd2}, and with time steps as indicated. The grid spacing is chosen as $h = 1.05 \sqrt{E/\rho} \, k$. }
    \label{fig:convString}
\end{figure}

In Figure \ref{fig:convString}, the output waveform of scheme \eqref{fd2} are checked against a reference solution obtained using the St{\"o}rmer-Verlet algorithm \eqref{eq:SVstring}, for the large input amplitude $\alpha = 300$, indicating that the output of the two schemes converges to a common solution.

\begin{figure}[hbt!]
    \centering
    \includegraphics[width = 0.9\linewidth]{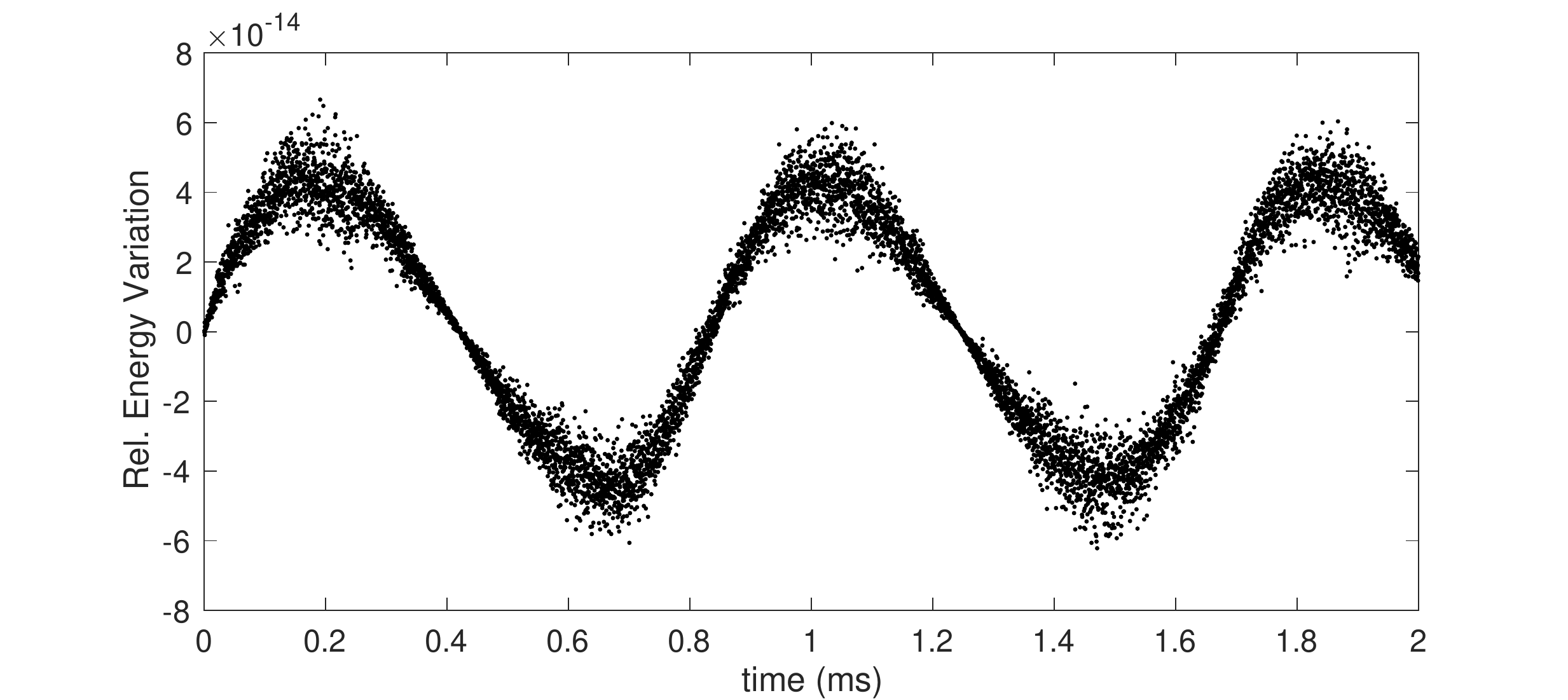}
    \caption{Energy variation of scheme \eqref{fd2} for the geometrically exact nonlinear string, using the split potential form as per \eqref{eq:EnStringSemiDiscSplit}. Here, the energy error is as per \eqref{relerrdef}. Here, the time step is $k=2.4 \cdot 10^{-7}$, and the grid spacing is chosen as $h = 1.05 \sqrt{E/\rho}\,k$. The string is initialised with using $\alpha = 300$. }
    \label{fig:stringEnErr}
 \end{figure}

The relative energy error, as defined in \eqref{relerrdef}, is shown in Figure \ref{fig:stringEnErr}, showing conservation to near machine accuracy. The larger range of variation here, compared with the case of the Fermi-Pasta-Ulam system is a result of the much larger state size and the resulting accumulation of errors. Notice in particular that here, in contrast to the case of the Fermi-Pasta-Ulam system, even though the energy variation is extremely small, there is now a clear correlation with the numerical solution. Such an effect is highly dependent on finite wordlength effects in double precision floating point, including the precise order of operations in the final update, and within the expression used to calculate the relative error \eqref{relerrdef}, and is not well understood. 

\begin{figure}[hbt!]
     \centering
     \includegraphics[width = 0.9\linewidth]{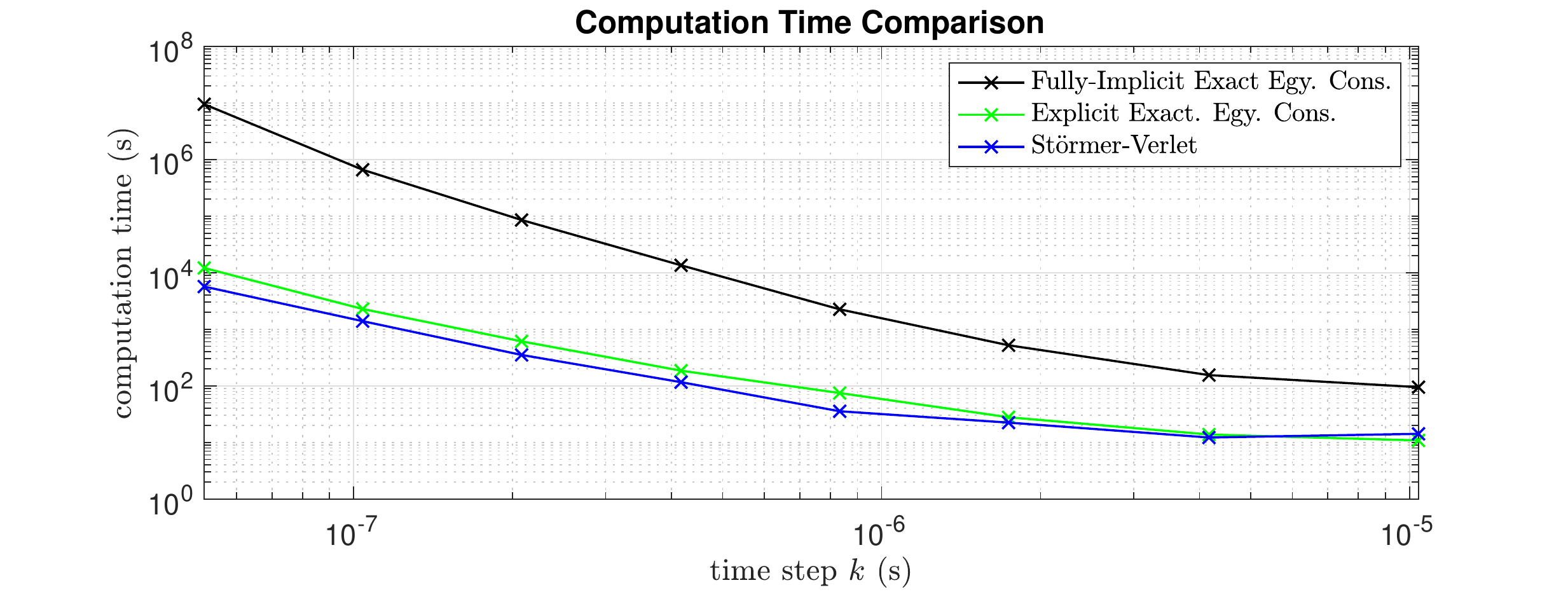}
     \caption{Computation times, in s, for the geometrically exact nonlinear string, using the fully implicit scheme \eqref{eq:NTstring} (black), St{\"o}rmer-Verlet \eqref{eq:SVstring} (blue) and the proposed scheme \eqref{fd2} (green). The grid spacing for a given time step is chosen as $h = 1.05\sqrt{E/\rho}\,k$. For the fully implicit scheme, Newton-Raphson is run with a tolerance of $10^{-13}$, and the maximum number of iterations is capped at 20.} 
     \label{fig:computeTimesString}
 \end{figure}

 Figure \ref{fig:computeTimesString} presents the compute times for schemes \eqref{eq:SVstring}, \eqref{eq:NTstring}  and \eqref{fd2}. The fully implicit, conservative scheme \eqref{eq:NTstring} is  the slowest, requiring a few iterations of the Newton-Raphson routine per time step. The proposed scheme \eqref{fd2} has compute times of the same order of magnitude as the fully explicit St{\"o}rmer-Verlet algorithm, and a few orders of magnitude smaller than the fully implicit scheme. All simulations were run in Matlab, using a 2016 MacBook Pro with a 2.9 GHz Quad-Core Intel Core i7 chip. 
 

\subsection{Nonlinear Plate Vibration: The F{\"o}ppl-von K{\'a}rm{\'a}n System}
\label{platesec}
As a final example, consider the problem of the transverse vibration of a thin plate at high amplitudes. A commonly used model is the so-called dynamic analogue of the F{\"o}ppl-von K{\'a}rm{\'a}n equations (see, e.g., \cite{vonkarman10, foppl1907, Nayfeh}), that have been used extensively recently in studies of wave turbulence \cite{ducceschi_PhysD_2014,yokoyama_PRL_2013,during2017wave}.
Time-stepping methods have been developed \cite{Kirby}, including a linearly-implicit energy-conserving method \cite{Bilbaonmpde08}. 

Though the dynamic F{\"o}ppl-von K{\'a}rm{\'a}n equations can be written directly in Hamiltonian form, they are most commonly presented as the following pair of coupled partial differential equations:
\begin{equation}
\label{vkdef}
    \rho\xi \partial_{t}^2 q  = -Q\Delta\Delta q +{\mathcal L}(q,F)\qquad\qquad \frac{2}{E\xi}\Delta\Delta F = -{\mathcal L}(q,q)\, .
\end{equation}
Here, $q(x,y,t)$ and $F(x,y,t)$ are the transverse displacement of the plate and Airy stress function respectively; both are functions of spatial coordinates $(x,y)\in{\mathcal D}\subset{\mathbb R}^{2}$, 
and time $t\geq 0$. In this simple example, the spatial domain ${\mathcal D}$ will be taken to be the square of side length $L$ m, so that ${\mathcal D} = [0,L]^2$. $\partial_{t}$ and $\Delta$ represent partial differentiation with respect to time $t$ and the 2D Laplacian operator, respectively. $\Delta\Delta$ is thus the biharmonic operator. For simplicity, boundary conditions are assumed to be of simply supported type over the boundary $\partial {\mathcal D}$ of ${\mathcal D}$, so that 
\begin{equation}
\label{ssdef}
    q = \Delta q = 0\qquad\qquad F = \Delta F = 0\qquad {\rm over} \quad\partial{\mathcal D}\, .
\end{equation}

The various constants that appear in \eqref{vkdef} are: $\rho$, the material density, in kg$\cdot$ m$^{-3}$; $\xi$, the plate thickness, in m; $E$, Young's modulus, in kg$\cdot$ s$^{-2}$ $\cdot$ m$^{-1}$; and the flexural rigidity $Q = E\xi^3/12(1-\nu^2)$, where $\nu$ is Poisson's ratio for the plate material. ${\mathcal L}$ is a bilinear operator, defined in terms of its action on two functions $f(x,y)$ and $g(x,y)$ as
\begin{equation}
\label{Ldef}
    {\mathcal L}(f,g) = \partial_{x}^2 f\partial_{y}^2 g+\partial_{y}^2 f\partial_{x}^2 g-2\partial_{x}\partial_{y}f\partial_{x}\partial_{y}g\, ,
\end{equation}
where $\partial_{x}$ and $\partial_{y}$ represent partial differentiation with respect to $x$ and $y$, respectively. Notice that only the first of \eqref{vkdef} is dynamic; the pair of equations \eqref{vkdef} could be rewritten as a single equation in displacement $q$ alone, and would constitute a second order in time cubic nonlinear PDE. 

System \eqref{vkdef} is Hamiltonian, with the total energy defined by
\begin{equation}
    H = \iint_{{\mathcal D}}\frac{\rho \xi}{2}(\partial_{t}q)^2 +\frac{Q}{2}(\Delta q)^2+\frac{1}{2E\xi}(\Delta F)^2d\sigma\, .
\end{equation}
This particular form of the energy holds under fixed edge boundary conditions, (such as simply supported \eqref{ssdef}), and must be modified under other types of conditions, such as free-edge. Notice that the final two terms under the integral above, which correspond to the potential energy, separate into a quadratic form in $q$, representing stored energy due to linear effects, and a quadratic form in $F$, representing additional nonlinear effects; both terms are individually non-negative, signalling that in numerical design, a splitting of the form described in Section \ref{splitting_sec} is available. 

\subsubsection{Semi-discrete Form}

For the square region ${\mathcal D}$, of side length $L$, one may start by defining grid locations $x_{l} = lh$, $y_{m} = mh$, where $l,m=1,\hdots,M-1$, for some integer $M$ such that $M = L/h$, where $h$ is a grid spacing. The semi-discrete grid functions $q_{l,m}(t)$ and $F_{l,m}(t)$ thus represent approximations to $q(x,y,t)$ and $F(x,y,t)$ at $x=x_{l}$ and $y=y_{m}$, respectively. 

Approximations to the spatial derivative operators $\partial_{x}$ and $\partial_{y}$ may be written, in terms of their action on a grid function $u_{l,m}(t)$ (such as $q_{l,m}$ or $F_{l,m}$ as defined above), as
\begin{equation}
\label{Dfirstdef}
    D_{x\pm}u_{l,m} = \mp\frac{1}{h}\left(u_{l,m}-u_{l\pm 1,m}\right)\qquad\qquad D_{y\pm}u_{l,m} = \mp\frac{1}{h}\left(u_{l,m}-u_{l,m\pm 1}\right)\, .
\end{equation}
Under simply supported conditions, when grid points outside the range $l,m=1,\hdots,M-1$ are referred to in the definitions above, such values are assumed to be zero. These are the most basic forward and backward approximations to derivatives, but are sufficient for the present purposes---more elaborate approximations (such as those of spectral type \cite{Vichnevetsky}) are available.

From the basic operations defined in \eqref{Dfirstdef}, centered approximations to the Laplacian $\Delta$ and biharmonic operator $\Delta\Delta$ follow as
\begin{equation}
    D_{\Delta} = D_{x+}D_{x-}+D_{y+}D_{y-}\qquad\qquad D_{\Delta\Delta} = D_{\Delta}D_{\Delta}\, .
\end{equation}
It is important to note that this particular construction of the biharmonic approximation $D_{\Delta\Delta}$, through a product of Laplacian approximations $D_{\Delta}$ under fixed conditions ensures that the simply supported conditions \eqref{ssdef} are satisfied. 

A semi-discrete approximation to \eqref{vkdef} then follows as
\begin{equation}
\label{vk2}
    \rho\xi\ddot{q}_{l,m} = -QD_{\Delta\Delta}q_{l,m}+\ell (q,F) = 0\qquad\qquad \frac{2}{E\xi}D_{\Delta\Delta}F_{l,m} = -\ell (q,q)\, .
\end{equation}
The operator $\ell(\cdot,\cdot)$ approximates ${\mathcal L}(\cdot,\cdot)$, as defined in \eqref{Ldef}. One useful centered approximation, operating on two grid functions $f_{l,m}$ and $g_{l,m}$ is \cite{Bilbaonmpde08}:
\begin{eqnarray*}
    \ell(f,g) &=& D_{x+}D_{x-}fD_{y+}D_{y-}g+D_{y+}D_{y-}fD_{x+}D_{x-}g\\
    && -\tfrac{1}{2}D_{x+}D_{y+}fD_{x+}D_{y+}g-\tfrac{1}{2}D_{x+}D_{y-}fD_{x+}D_{y-}g-\\ &&\tfrac{1}{2}D_{x-}D_{y+}fD_{x-}D_{y+}g-\tfrac{1}{2}D_{x-}D_{y-}fD_{x-}D_{y-}g \, .\notag
\end{eqnarray*}

It is useful to represent the semi-discrete ODE system \eqref{vk2} in vector form, using $(M-1)^2\times 1$ vectors ${\bf q}$ and ${\bf F}$. In this representation, the operators $D_{x+}$ and $D_{y+}$ become $(M-1)^2\times M(M-1)$ matrices ${\bf D}_{x+}$ and ${\bf D}_{y+}$, and ${\bf D}_{x-} = -{\bf D}_{x+}^{\intercal}$ and ${\bf D}_{y-} = -{\bf D}_{y+}^{\intercal}$. The operators $D_{\Delta}$ and $D_{\Delta\Delta}$ become $(M-1)^2\times (M-1)^2$ matrices ${\bf D}_{\Delta}$ and ${\bf D}_{\Delta\Delta}$ respectively. One arrives at the form 
\begin{equation}
\label{vk3}
    \rho\xi\ddot{{\bf q}}= -Q{\bf D}_{\Delta\Delta}{\bf q}+\ell ({\bf q},{\bf F}) = 0\qquad\qquad \frac{2}{E\xi}{\bf D}_{\Delta\Delta}{\bf F} = -\ell ({\bf q},{\bf q})\, .
\end{equation}
This second order in time ODE system serves as the starting point for methods such as St{\"o}rmer-Verlet and other linearly-implicit energy-conserving methods, as described below. 

By introducing the momentum variable ${\bf p} = {\bf M}\dot{{\bf q}}$, the first order system \eqref{Hamdyn} equivalent to \eqref{vk3} results from a  Hamiltonian of the form of \eqref{Hamdef}, with $N=(M-1)^2$, and where
\begin{equation}
\label{vk_VMdef}
    {\bf M} = \rho\xi h^2{\bf I}_{(M-1)^2}\qquad V = \frac{Qh^2}{2}\|{\bf D}_{\Delta}{\bf q}\|^2+\frac{h^2}{2E\xi}\|{\bf D}_{\Delta}{\bf F}\|^2\, .
\end{equation}
Here, ${\bf I}_{(M-1)^2}$ is the $(M-1)^2\times (M-1)^2$ identity matrix. A natural splitting of the potential energy $V$ as in \eqref{split_def} follows, with 
\begin{equation}
\label{vksplit}
    {\bf K} = Qh^2{\bf D}_{\Delta\Delta}\qquad\qquad V' = \frac{h^2}{2E\xi}\|{\bf D}_{\Delta}{\bf F}\|^2\, .
\end{equation}

\subsubsection{Numerical Methods}

Beginning from the second order system \eqref{vk3}, one may introduce the discrete time vectors ${\bf q}^{n}$ and ${\bf F}^{n}$, representing approximations to $q$ and $F$ at $t=nk$, for integer $n$, and where $k$ is the time step in s. St{\"o}rmer-Verlet integration results immediately in:
\begin{equation}
\label{vk_stormer}
    {\bf q}^{n+1} = \left(2{\bf I}_{(M-1)^2}-\frac{Qk^2}{\rho\xi}{\bf D}_{\Delta\Delta}\right){\bf q}^{n}-{\bf q}^{n-1}  +\frac{k^2}{\rho\xi}\ell ({\bf q}^{n}, {\bf F}^{n}) \qquad\qquad \frac{2}{E\xi}{\bf D}_{\Delta\Delta}{\bf F}^{n} = -\ell ({\bf q}^{n}, {\bf q}^{n})\, .
\end{equation}
The first of these updates is explicit, but relies on ${\bf F}^{n}$, which must be obtained from the second equation through the solution of a linear system involving the biharmonic operator ${\bf D}_{\Delta\Delta}$.

The St{\"o}rmer-Verlet scheme above is not conservative, and is prone to numerical instability. A slight variant, however, leads to a scheme with exact energy conservation:
\begin{subequations}\label{vk_oldegy}
\begin{align}
    &{\bf q}^{n+1} = \left(2{\bf I}_{(M-1)^2}-\frac{Qk^2}{\rho\xi}{\bf D}_{\Delta\Delta}\right){\bf q}^{n}-{\bf q}^{n-1}  +\frac{k^2}{2\rho\xi}\ell ({\bf q}^{n},{\bf F}^{n+1}+{\bf F}^{n-1}) \\ &\frac{1}{E\xi}{\bf D}_{\Delta\Delta}\left({\bf F}^{n+1}+{\bf F}^{n}\right) = -\ell ({\bf q}^{n+1}, {\bf q}^{n})\, .
\end{align}
\end{subequations}
Due to the bilinearity of the operator $\ell$, this scheme is linearly implicit---at each time step, ${\bf q}^{n+1}$ and ${\bf F}^{n+1}$ must be solved simultaneously, using a linear system constructed anew at each time step; St{\"o}rmer-Verlet, in contrast, requires only the solution of a linear system in ${\bf D}_{\Delta\Delta}$, which is of a known form. Thus it may be expected that scheme \eqref{vk_oldegy}, while energy-conserving and provably numerically stable \cite{Bilbaonmpde08}, is significantly slower to execute than the St{\"o}rmer-Verlet scheme \eqref{vk_stormer}. 

Finally, from the Hamiltonian form given in \eqref{Hamdyn}, with $V$ and ${\bf M}$ as given in \eqref{vk_VMdef}, as well as the splitting of the potential energy as in \eqref{vksplit}, a time-interleaved scheme of the form \eqref{fd2} results. This scheme possesses a non-negative exactly conserved numerical energy under the condition
\eqref{stabcond} which, in this case, reduces to a lower bound on the grid spacing $h$ in terms of the time step $k$:
\begin{equation}
\label{vk_stabcond}
    h\geq h_{{\rm min}} = 2\sqrt{k}\left(D/\rho \xi\right)^{\tfrac{1}{4}}\, .
\end{equation}
The scheme, like St{\"o}rmer-Verlet, relies on a linear system solution involving the biharmonic operator ${\bf D}_{\Delta\Delta}$, but is otherwise explicit. 

\subsubsection{Numerical Examples}

As an example, consider initialisation of the F{\"o}ppl-von K{\'a}rm{\'a}n system \eqref{vkdef} through its lowest linear mode shape:
\begin{equation}
\label{vk_ic}
    q(x,y,0) = \alpha\xi\sin(\pi x/L)\sin(\pi y/L)\qquad \partial_{t}q|_{x,y,t=0} = 0\, ,
\end{equation}
where the dimensionless parameter $\alpha$ controls the maximum amplitude of the initial condition relative to the plate thickness $\xi$. Furthermore, the plate is assumed to be made of steel with $E=2\times 10^{11}$ Pa, $\rho = 7850$ kg$\cdot$ m$^{-3}$, and $\nu = 0.3$, and to be of thickness $\xi = 2$ mm and side length $L=0.5$ m. 

\begin{figure}[hbt!]
\centering
\includegraphics[width = 0.9\linewidth]{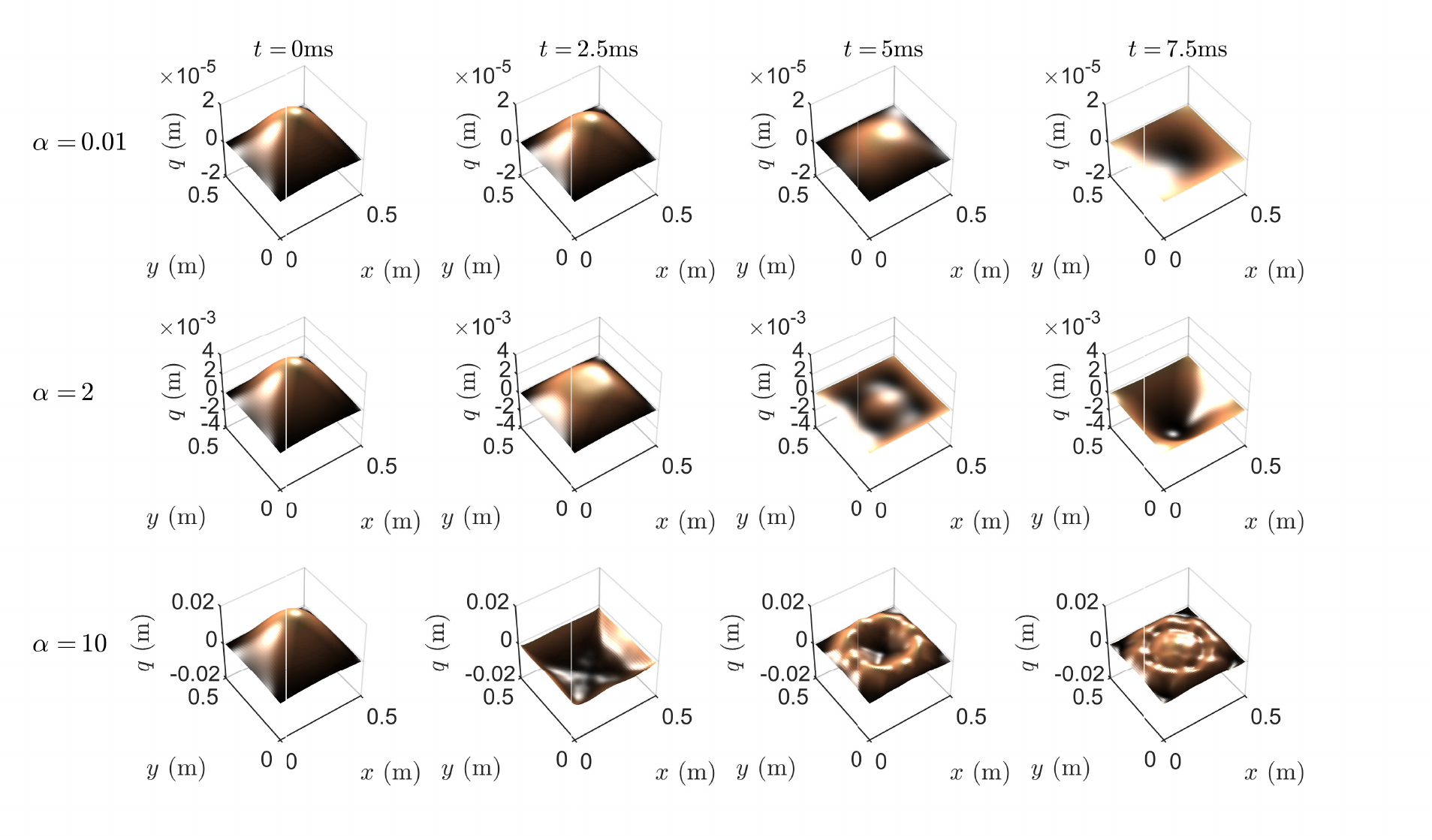}
\caption{Snapshots of the time evolution of a F{\"o}ppl-von K{\'a}rm{\'a}n plate, at times as indicated, and for different initial condition amplitudes $\alpha = 0.01$, $\alpha = 2$ and $\alpha = 10$. Results are computed using scheme \eqref{fd2}, with a time step of $k=5\times 10^{-6}$. }
\label{vkevfig}
\end{figure}

\begin{figure}[hbt!]
\centering
\includegraphics[width = 0.9\linewidth]{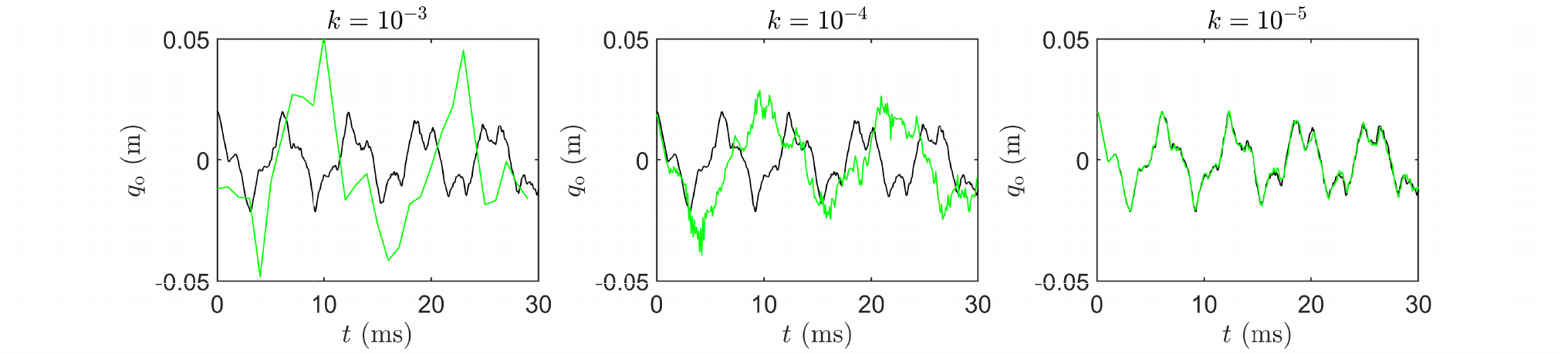}
\caption{The output waveform $q_{{\rm o}}(t)$ drawn from a plate vibrating at high amplitude, with $\alpha = 10$. The reference solution, generated using St{\"o}rmer-Verlet with $k=2.5\times 10^{-6}$ s is shown in black; results are shown for scheme \eqref{fd2}, in green, at different time steps $k$, as indicated. }
\label{vkconvfig}
\end{figure}

Using the scheme \eqref{fd2}, with the grid spacing and time step chosen according to \eqref{vk_stabcond}, typical amplitude-dependent  behaviour is observed. See Figure \ref{vkevfig}. At a low initial condition amplitude of $\alpha = 0.01$, behaviour is essentially linear. At amplitudes near the plate thickness at $\alpha = 2$, the period of oscillation decreases, and spontaneous mode generation is observed, and for large amplitudes, such as $\alpha = 10$, turbulent behaviour is observed. Under such stringent high-amplitude conditions ($\alpha = 10$), scheme \eqref{fd2} remains stable, and results converge to those of the reference solution, computed using St{\"o}rmer-Verlet with a time step of $k=2.5\times 10^{-6}$ s. No regularisation (see Section \ref{shift_sec}) is used in this case. See Figure \ref{vkconvfig}. Under these conditions, St{\"o}rmer-Verlet is unstable for $k>2\times 10^{-5}$ s. 

\begin{figure}[hbt!]
\centering
\includegraphics[width = 0.9\linewidth]{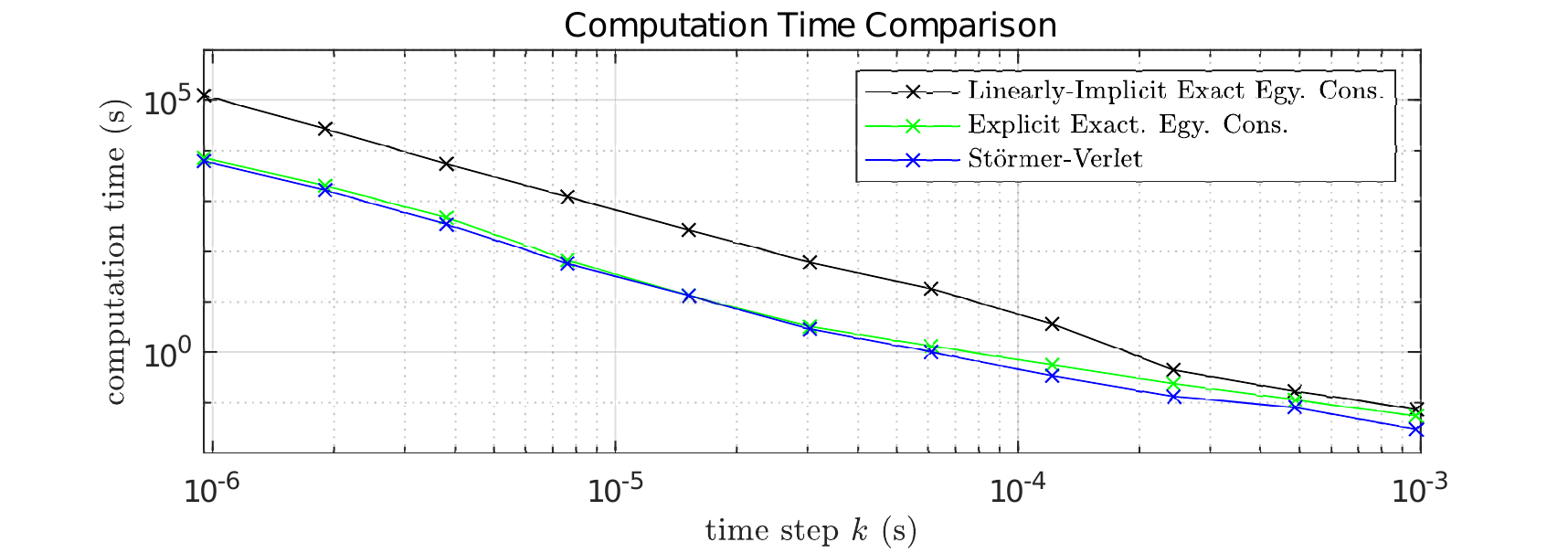}
\caption{Computation times, in s, for the  F{\"o}ppl-von K{\'a}rm{\'a}n plate, for a 1 s simulation duration, using the linearly implicit exact energy conserving method \eqref{vk_oldegy}, the explicit exact energy conserving method \eqref{fd2}, and St{\"o}rmer-Verlet \eqref{vk_stormer}.  }
\label{vktimesfig}
\end{figure}

The simulation of the F{\"o}ppl-von K{\'a}rm{\'a}n system is computationally intensive. Most interesting in this case is a comparison of computation times between St{\"o}rmer-Verlet \eqref{vk_stormer}, the linearly-implicit energy-conserving method \eqref{vk_oldegy}, and the explicit energy-conserving scheme \eqref{fd2}. Unavoidable in all cases is some form of linear system solution; for \eqref{vk_stormer} and \eqref{fd2}, this will involve the biharmonic operator, which is constant over the course of a simulation, and thus amenable to factorisation techniques (e.g. Cholesky) to decrease solution times. For the linearly-implicit exact energy conserving method \eqref{vk_oldegy}, however, this is not the case---the linear system to be solved must be constructed anew at each time step. This is reflected in timings, as shown in Figure \ref{vktimesfig}, for simulations of 1 second for different choices of time step $k$.  Computation was performed in Matlab on a Lenovo P50 with an Intel Xeon E3 v5. As can be seen, computation time for the explicit scheme \eqref{fd2} is far lower than for scheme \eqref{vk_oldegy}, and very nearly on par with St{\"o}rmer-Verlet. 

\begin{figure}[hbt!]
\centering
\includegraphics[width = 0.9\linewidth]{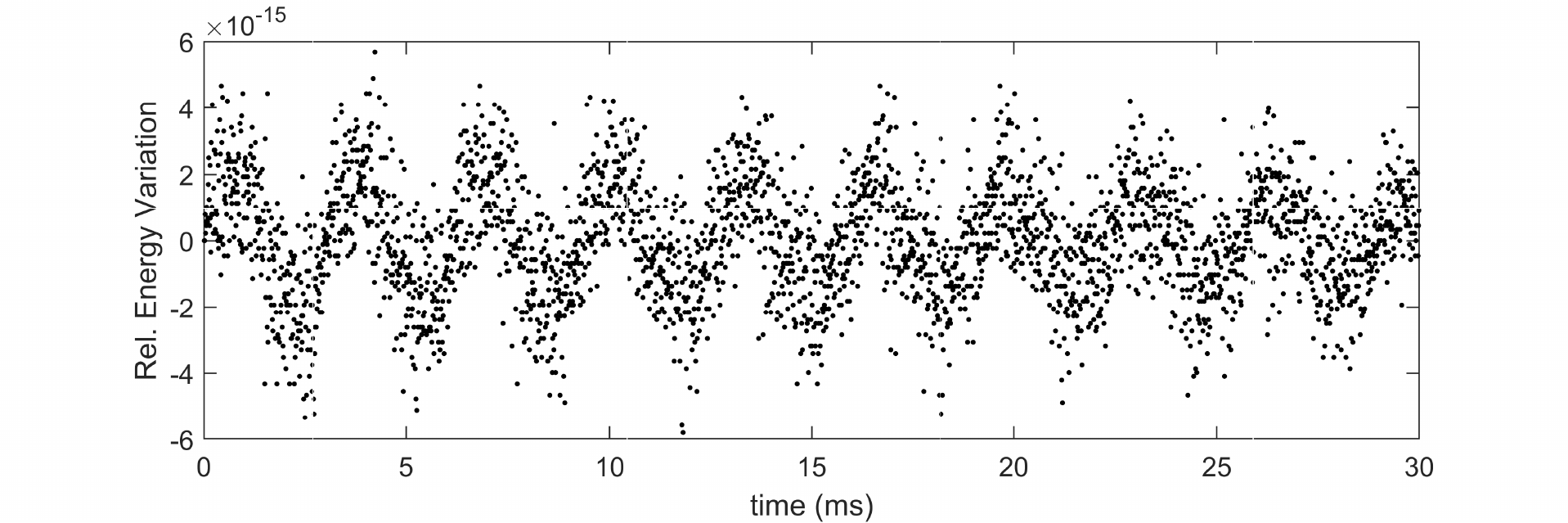}
\caption{Relative numerical energy variation, as defined in \eqref{relerrdef}, for scheme \eqref{fd2} for the F{\"o}ppl-von K{\'a}rm{\'a}n plate, with a high initial condition of the form of \eqref{vk_ic}, with $\alpha = 10$. The time step is $k= 10^{-5}$ s. }
\label{vkegyfig}
\end{figure}

Finally, see Figure \ref{vkegyfig}, illustrating the relative energy variation for scheme \eqref{fd2} for the F{\"o}ppl-von K{\'a}rm{\'a}n plate; as in the case of nonlinear string vibration, the energy variation is of the order of $10^{-15}$, with some correlation with the numerical solution visible---see Figure \ref{fig:stringEnErr} for comparison. 

\section{Concluding Remarks}
\label{concsec}

The design of exactly energy-conserving methods for Hamiltonian systems has progressed from fully implicit designs through, more recently, to explicit methods for which exact energy conservation can be attained through an approximation to a continuous integral of the potential energy, or, more importantly, to linearly implicit designs based on invariant energy quadratisation. The main novelty in this paper is to call attention to structure within such linearly implicit designs---structure that can be exploited in order arrive at fully explicit methods. These exactly conservative methods are of roughly the same computational cost as the most efficient non-conservative explicit methods---with the additional feature of a clear means of ensuring numerical stability, either unconditionally, or, if a splitting of the potential energy is employed, under well-defined conditions on the time step that are independent of the initial conditions. Accuracy is of second order, and is plainly evident in the centered (but interleaved) discretisation approach, and borne out by simulation results (see Section \ref{fpusec}). It is not clear whether it is possible to extend this framework to obtain higher order accuracy.  

These methods are not completely general, and require, additionally, a condition of non-negativity on the potential energy, as per invariant energy quadratisation approaches. More generally, given that the dynamics of a system are independent of shifts in the potential energy by a constant (i.e., a gauge), a more general condition is that the potential energy is bounded from below. The useful technique of splitting of the potential energy introduced here is a further restriction. But the restriction \eqref{Vpos} mentioned above to non-negative expressions for potential energy $V$ (or bounded from below) is slightly more strict than necessary. More general is a restriction to expressions $V({\bf q})$ that are single-signed for all ${\bf q}$---and even more generally, bounded either from above or below. An important example here is the $N$-body problem, under a gravitational potential,  for which $V({\bf q})\leq 0$. In this case, one may set, instead of \eqref{psidef}, $V = -\tfrac{1}{2}\psi^2$, and the main development follows as above, with this sign change, and an explicit exactly energy conserving method follows as before. In this case, however, global bounds on solution size are not available, as the total energy itself is no longer necessarily non-negative. See the comments in \cite{Gonzalez96} regarding general dynamic stability for Hamiltonian systems.  

Only briefly alluded to here, in Section \ref{remarkssec}, is the extension to the case of variable time steps---useful in modeling systems with slow/fast dynamics, as discussed in \cite{Marazzato}. The behaviour of such schemes remains unexplored. Another open question is the need for regularization, as introduced in Section \ref{shift_sec}. In two of the three cases presented here, the Fermi-Pasta-Ulam problem, and the F{\"o}ppl-von K{\'a}rm{\'a}n system, good results were obtained without the use of such regularisation. In the case of string vibration, however, such regularisation was necessary. It would be of great interest to know the origin of this distinction. In all cases, however, regularisation has no impact on the explicit nature of the schemes proposed here, and only negligible impact on computational cost. 

\section*{Acknowledgments}
M. Ducceschi was supported by the European Research Council (ERC), under grant 2020-StG-950084-NEMUS. For the purpose of open access, the first author has applied a creative commons attribution (CC BY) licence to any author accepted manuscript version arising.

\bibliographystyle{unsrt}
\bibliography{bilbao_submission_arxiv}

\end{document}